\documentclass[11pt]{conm-p-l}

\usepackage{epsf,amsmath,amssymb,amscd,amsthm,color,fancyhdr}

\usepackage{times}
 \usepackage[all] {xy}

\usepackage{pdfsync}

\newtheorem{theo}{Theorem}[section]

\newtheorem{prop}[theo]{Proposition}

\newtheorem{cor}[theo]{Corollary}

\newtheorem{defi}[theo]{Definition}

\newtheorem{rema}[theo]{Remark}

\newcommand{\qeed}{\hfill\textrm{QED}\break\null}
\newenvironment{demo}{\noindent\textit{Proof.}~}{\qeed}

\begin{document}

\title[Pure Spinor Algebra]{ Pure spinors and a construction of the $E_*$- Lie algebras }

\author[Slupinski]{Marcus J. Slupinski}
\address{IRMA, Universit\'e de Strasbourg\\
 7 rue Ren\'e Descartes\\
F-67084 Strasbourg Cedex FRANCE}
\email{marcus.slupinski@math.unistra.fr} 
\thanks{The first author acknowledges the recurring support from the Math. Research Inst., OSU, of our collaboration}
 
\author[Stanton]{Robert J. Stanton}
\address{1937 Beverly Rd, Columbus, OH 43221}
\email{stanton.2@osu.edu}
\thanks{The second author is grateful to IRMA, Universit\'e de Strasbourg, for its hospitable and stimulating environment.}


\dedicatory{This paper is dedicated to our friend and collaborator Gestur \'Olafsson.}


\begin{abstract}
Let $(V,g)$ be a $2n$-dimensional hyperbolic space and $C(V,g)$ its Clifford algebra.  $C(V,g)$ has a $\mathbb Z$-grading, $C^k $, and an algebra isomorphism $C(V,g)\cong End(S)$, $S$ the space of spinors. \'E. Cartan defined operators $L_k: End(S) \to C^k$ which are involved in the definition of pure spinors. We shall give a more refined study of the operator $L_2$, in fact, obtain explicit formulae for it in terms of spinor inner products and combinatorics, as well as the matrix of it in a basis of pure spinors. Using this information we give a construction of the exceptional Lie algebras $\mathfrak e_6, \mathfrak e_7, \mathfrak e_8$  completely within the theory of Clifford algebras and spinors.
 \end{abstract}

\date{\today}

\maketitle
       
\section{Introduction}
Constructions of exceptional Lie algebras over quite general fields have been given by many people and from various perspectives. While the list is too long to give, error free, we must mention Freudenthal and Tits. The perspective of this paper is that of spinors. When the base field $k$ is $\mathbb R$, the classification by \'E. Cartan of irreducible Riemannian symmetric spaces of the noncompact type already provides an example, the case of $\mathfrak e_8$, for which the pair ($\mathfrak {so}(16,\mathbb R), S_{\pm})$ occurs. More recently, J. F. Adams \cite {Ad} gave a construction of compact exceptional Lie groups using compact spin groups and the relationship of some to Jordan algebras. Moroianu-Semmelman \cite {MoSe}  gave a construction of exceptional Lie algebras of compact type by refining Kostant's \cite{K} invariant 4-tensor characterization of certain holonomy representations and coupled with the compact spinor material from \cite {Ad}. Our point of view is to present natural properties of Clifford algebras and their spinors for a hyperbolic space over a very general field $k$, and then to derive the existence of the $E_*$-series using these properties and the combinatorics of pure spinors, thus a construction intrinsic to spinor algebra. 

The paper is essentially self-contained and written with Lie theorists in mind, such as a master like Gestur, hence includes some standard material on spinors known to experts. We begin with a Clifford algebra $C(V,g)$ and its basic isomorphism with $End(S)$, $S$ the spinors. Then we relate properties of $C(V,g)$ and $S$, including fundamental material about the spinor norm. The key tool in the paper is the operator introduced by Cartan, $L_2: S\times S \to C_2$. After we give a new description of $L_2$,  we obtain an explicit formula for $L_2$ in terms of Clifford elements. The formula is intrinsic to Clifford theory as the operator is completely specified by various spinor norms and combinatorics. Then using the basis of pure spinors we compute the matrix of $L_2$ and express its entries in terms of spinor norms and combinatorics. This treatment is completely general for a hyperbolic space $(V,g)$ and a field $k$ of characteristic not $2$ or $3$. 

In the last section we specialise the formula for the matrix of $L_2$ to three specific dimensions and show that various entries of the matrix vanish for combinatorial reasons yielding a Jacobi identity for the various Lie algebras in the $E_*$-series.

There are several potential future directions. The choice of a hyperbolic $g$ was made to avoid field extensions of $k$ - indeed there are metrics of other signature that could be considered. Also, the combinatorics that arise in the computations mirror properties of the Weyl group quotient that parametrises Schubert cells in the flag variety of projectivised pure spinors. We did not consider whether other topological properties of the cells are responsible for the various combinatorial identities. Finally, since the spinor algebra is a universal linear construction, we expect the spinor algebra constructions, in particular $L_2$, to transfer to vector bundles.

\section{Background on Spinors } 
Let $V$ be a $2n$-dimensional  vector space over a field $k$ of characteristic not $2$ or $3$. We shall assume that $V$ has a nondegenerate symmetric bilinear form $g$ of Witt index $n$ (i.e. a hyperbolic form). The hyperbolic case allows us to give a rather complete presentation of the results without any base extension of $k$. This was highlighted by Chevalley and today seems even more relevant. A good reference for much of the basic material of this section is \cite{Ch}.  
\subsection{Clifford algebra}\hfill
\vskip 0.1in

Let $C = C(V,g) $ be the Clifford algebra of $V$ with respect to  $g$. Then $C$  has the usual ${\mathbb Z}_2$ grading $C=C_+\oplus C_-$ inherited from the tensor algebra of $V$. 
 As $g$ is hyperbolic,  $C$ is isomorphic to the algebra of
$2^n\times2^n$  matrices over $k$. We can choose a $2^n$-dimensional $k$-vector space $S$, up to equivalence, called the space of spinors,  and obtain an
algebra isomorphism
$$
C\cong End(S).
$$
Hence $C$ has a natural trace that we denote $Tr:C(V,g) \rightarrow {k}$.

The vector space $V$ is naturally included in $C_-$ so,  from now on, we  consider $V$ as a subset
of $C$. By the universal property of  $\Lambda^* (V)$, the exterior algebra of $V$,  one can extend the inclusion of $V$ into $C_-$ to an $O(V,g)$- equivariant linear (but not algebra) isomorphism  $Q:\Lambda^* (V)
\rightarrow C$ by defining
$$
Q(v_1\wedge\dots\wedge v_k)=\frac{1}{k!}\sum_{\sigma}(-1)^{\sigma}v_{\sigma(1)}\cdots v_{\sigma(k)},
$$
at least if $k$ is of characteristic $0$.
If $v_{1},\dots ,v_{k}$ are orthogonal  this formula implies that
$$
Q(v_{1}\wedge\dots\wedge v_{k})=v_{1}\cdots v_{k}
$$
and one uses  this property to  characterize $Q$ when $k$ is of positive characteristic (see \cite{Ch}).

We set $C^k = Q(\Lambda^k (V))$. Hence if $\{v_{1},\dots ,v_{2n}\}$ is a basis of $V$ and $i_1<i_2<\cdots < i_k$ then the collection $\{Q(v_{i_1}\wedge\dots\wedge v_{i_k})
\}$ is a basis of $C^k$. The collection of subspaces $C^k$ then give $C$ the structure of a $\mathbb Z$ graded vector space.

 $C$ is also a filtered algebra where $D^k$ is generated by products of at most $k$ elements of $V$. We have then an isomorphism of the associated $\mathbb Z$ graded space determined by the filtration onto the $\mathbb Z$ graded  $Q(\Lambda^* (V))$. 

The following commutator relations are well known:
$$
[C^1,C^1]\subseteq C^2,\quad[C^2,C^m]\subseteq C^{m},\quad [C^2,C^{2n}]=0.
$$
Consequently, $C^2$, $C^1\oplus C^2$ and $C^2\oplus C^{2n}$ are Lie algebras. 

The composition with $Q$ of any
$\mathfrak {o}(V,g)$-equivariant isomorphism
$$
\mathfrak {o}(V,g)\cong \Lambda^2(V)
$$
defines a Lie algebra  isomorphism $\mathfrak {o}(V,g)\cong C^2$. 

Similarly one shows easily that the orthogonal Lie algebra of a vector space of dimension $2n+1$  of maximal Witt index is isomorphic to $C^1\oplus C^2$.

The canonical anti-automorphism of order 2 of $C(V,g)$, namely the one extending $v\mapsto v$ for $v\in V$, is inherited from the tensor algebra. It will be denoted $x\mapsto x^T$. Using the canonical anti-automorphism $T$ and the trace $Tr$ one can give $C$ a norm, namely $\Vert c \Vert^2 = Tr(c^Tc)$.

\begin{prop}
 Let $g_{\Lambda}$ be the natural extension of $g$ to $\Lambda^*(V)$.
Then
$$
2^ng_{\Lambda}(\alpha,\beta)=Tr(Q(\alpha)^TQ(\beta)),
$$
i.e. $Q$ is a multiple of an isometry.

\end{prop}
\begin{rema}
By an orthonormal basis of $V$ we mean a basis $\{e_1,\dots ,e_{2n}\}$  which satisfies 
$$
g(e_i,e_j)=\pm \delta_{ij}.$$
Orthonormal bases exist because $(V,g)$ is isometric to an orthogonal  sum of hyperbolic planes. Later we will use ordered orthonormal bases.
\end{rema}
For the $\mathbb Z$ grading
$C=\oplus C^k$ and with respect to an orthonormal basis there is a  formula for the projection
$\pi_k:C\rightarrow C^k$: 
$$
\pi_k(c)=\frac{1}{2^{n}}\sum_{i_1<\dots <i_k}g(e_{i_1},e_{i_1})\dots g(e_{i_k},e_{i_k})Tr(e_{i_k}\dots e_{i_1}c)e_{i_1}\dots
e_{i_k},$$
as follows easily from the fact that $\{Q(v_{i_1}\wedge\dots\wedge v_{i_k})\}$ is a basis of $C^k$ and $Q(e_{i_1}\wedge\dots\wedge e_{i_k})=e_{i_1}\cdots e_{i_k}$.

There is a natural $\mathbb Z_2$ grading  of $S$, $S=S_1\oplus S_2$, into a direct sum of two $2^{n-1}$ dimensional subspaces compatible with the $\mathbb Z_2$ graded action  of $C$, i.e.,
$$
 C_+\cdot S_1\subseteq S_1,\quad C_+\cdot S_2\subseteq S_2,\quad C_-\cdot S_1\subseteq S_2\hbox{ and }C_-\cdot S_2\subseteq S_1.
$$
Elements of $S_1$ or $S_2$ are called half-spinors. An element of $C$  that implements such a grading of $S$, i.e. which is the identity on one
half-spinor space and minus the identity on the other, will be called a grading element. They are usually denoted by $\varepsilon$ and are elements of $C^{2n}$. Note that
if $\{e_1,\dots ,e_{2n}\}$ is an orthonormal basis of $V$ then
$$
\varepsilon=e_1\dots e_{2n}
$$
satisfies $\varepsilon^2=1$ and is a grading element. The graded decomposition of $S$ corresponding to this particular $\varepsilon$ will be denoted $S=S_+\oplus S_-$.

\begin{prop}
 If  $\varepsilon$ is a grading element, then for all $1\le k\le 2n$ and all $c\in C$,
 $$
 \pi_{2n-k}(\varepsilon c)=\varepsilon\pi_k(c).
 $$
 \end{prop}
 \begin{demo} This is essentially III.4.3 in \cite{Ch}.
 \end{demo}
\subsection{Spinor norms}\hfill 

Recall that $V$ is of dimension $2n$ and that $g$ is hyperbolic.

\begin{defi} A spinor norm is a  bilinear map 
$ B:S\times S\rightarrow k$  such that
$$
B(v\cdot\phi,\psi)=B(\phi,v\cdot\psi)\quad\forall v\in V,\forall \phi,\psi\in S .
$$
Hence, given a non-zero spinor norm $B$, the canonical anti-automorphism of $C$, $x\mapsto x^T$, is, via the isomorphism $C\cong End(S)$, the transpose relative to $B$.
\end{defi}

\begin{prop}(\'E. Cartan)
 The space of spinor norms is one dimensional.
\end{prop}

The  symmetry type of a spinor norm, symmetric or alternating, as well as the type of its restriction to the half-spinor spaces is summarized in the next result.

\begin{prop}\label{spinornormsymmetries}
 Let $B$ be a spinor norm.

\noindent(i) If $n\equiv 0$ mod $4$ then $B$ is symmetric and even, i.e., $B(S_1,S_2)=0$.

\noindent(ii) If $n\equiv 1$ mod $4$ then $B$ is symmetric and odd, i.e., $B(S_1,S_1)=B(S_2,S_2)=0$.

\noindent(iii) If $n\equiv 2$ mod $4$ then $B$ is antisymmetric  and even.

\noindent(iv) If $n\equiv 3$ mod $4$ then $B$ is antisymmetric and odd.
\end{prop}
We shall fix a non-degenerate spinor norm, $B$. It is natural to have the relationship of the $\mathbb Z$ grading of $C$ to $B$. 
If $\phi,\psi\in S$, $1\le k\le 2n$ and $v_1,\dots ,v_k\in V$ are orthogonal, it is clear that
$$
B(v_1\cdots v_k\cdot\phi,\psi)=(-1)^{\frac{1}{2}k(k-1)}B(\phi,v_1\cdots v_k\cdot\psi).
$$
Hence we have

\begin{cor}\label{spinornorminvprops}
 If $k\equiv 2$ or $3$ (mod 4), then 
$$
B(c\cdot\phi,\psi)+B(\phi,c\cdot\psi)=0\quad\forall c\in C^k,
$$
i.e.,  a spinor norm is invariant under the action of $C^k$ iff $k\equiv 2$ or $k\equiv 3$ (mod 4). 

\end{cor}
Since $B$ can have either symmetry type, symmetric or alternating, we denote by $\mathfrak{aut}(S,B)$ the endomorphisms of $S$ that leave invariant the spinor norm.  

\begin{cor} 
Identifying $C(V,g)$ with ${End}(S)$, we have
$$
\bigoplus_{k\equiv 2\text{ or }3\text{ (mod }4)}^{}\hskip-,6cm C^k= \mathfrak {aut}(S,B).
$$
\end{cor}
\begin{demo} By Corollary \ref{spinornorminvprops} the LHS is included in the RHS. The result  follows from a dimension count for the corresponding symmetry type of $B$:
$$
\sum_{k\equiv 2\text{ or }3\text{ (mod }4)}^{}\binom{2n}{k}=2^{n-1}(2^{n}-1)={\rm dim}\,\mathfrak {so}(S,B)\quad (n\equiv 0,1\mod4)
$$ 
and
$$
\sum_{k\equiv 2\text{ or }3\text{ (mod }4)}^{}\binom{2n}{k}=2^{n-1}(2^{n}+1)={\rm dim}\,\mathfrak{sp}(S,B)\quad (n\equiv 2,3\mod4).
$$ 
\end{demo}
\begin{cor}
 If $\varepsilon$ is any grading element and we set 
$$
\mathfrak{aut}_{\pm}(S,B)=\{c\in\mathfrak{aut}(S,B):\,c\varepsilon=\pm \varepsilon c\}
$$
then
$$
\bigoplus_{k\equiv 2\text{ (mod }4)}^{} \hskip-,3cm C^k= \mathfrak{aut}_+(S,B),\qquad
\bigoplus_{k\equiv 3\text{ (mod }4)}^{}\hskip-,3cm C^k= \mathfrak{aut}_-(S,B).
$$
\end{cor}
\begin{demo} This is immediate from the previous Corollary.
\end{demo}

From Prop. \ref{spinornormsymmetries} we see that the spinor norm $B$ is even if and only if $n$  is even. In this case $B$ restricts to nondegenerate forms $B_1$ and $B_2$ on the half-spinor spaces $S_1$ and $S_2$ respectively.  Thus for $n$ even and grading element $\varepsilon$, the Lie algebra $ \mathfrak{aut}_+(S,B)$ is a  direct product 
$$
\mathfrak{aut}_+(S,B)\cong\mathfrak{ aut}(S_1,B_1)\oplus \mathfrak{ aut}(S_2,B_2).
$$
To realise this decomposition of $\mathfrak{aut}_+(S,B)$ in the Clifford algebra we  use the grading element.
\begin{cor}
Let $n$ be even. Let $\varepsilon\in C(V,g)$ be a grading element and $S=S_1\oplus S_2$ the associated grading (the $\pm 1$ eigenspaces of $\varepsilon$). Set $\varepsilon_{\pm} =\frac{(1\pm\varepsilon)}{2}$.

(i) Then
$$
\bigoplus_{k\equiv 2\text{ (mod }4)}^{} \hskip-,3cm C^k=\left(\bigoplus_{k\equiv 2\text{ (mod }4),k\le n}^{}\hskip-,1cm C^k \varepsilon_+\right)\oplus
\left(\bigoplus_{k\equiv 2\text{ (mod }4),k\le n}^{}\hskip-,1cm C^k \varepsilon_-\right).
$$
(ii) With respect to the decomposition $S=S_1\oplus S_2$, the
two summands of (i)  are:
$$
\bigoplus_{k\equiv 2\text{ (mod }4),\, k\le n}\hskip-,1cm C^k \varepsilon_+=
\left(
\begin{array}{c|c}
\mathfrak{aut}(S_1,B_1) & \quad0 \quad\\ \hline
0 & \quad0 \quad \\
\end{array}
\right)
$$
and
$$
\bigoplus_{k\equiv 2\text{ (mod }4),\,k\le n}\hskip-,1cm C^k \varepsilon_-=
\left(
\begin{array}{c|c}
\quad 0 \quad&0\\ \hline
\quad 0 \quad& \mathfrak{aut}(S_2,B_2)  \\
\end{array}
\right)
$$
\end{cor}

\begin{demo}
Part (i) basically follows from the fact that
$$
1=\frac{1}{2}(1+\varepsilon)+\frac{1}{2}(1-\varepsilon)
$$
decomposes the identity of $C$ as a sum of two orthogonal idempotents. Part (ii) is straightforward.
\end{demo}
\begin{rema} 
 If $k<n$,  both $C^k \varepsilon_+$ and $C^k \varepsilon_-$ are isomorphic to $\Lambda^k(V)$ as absolutely  irreducible $\mathfrak{so}(V,g)$-representations. If $k=n$, \'E. Cartan showed that $C^n \varepsilon_+$ and $C^n \varepsilon_-$ are absolutely  irreducible, non-isomorphic representations of the same dimension. The proposition therefore gives an explicit reduction of $\mathfrak{aut}(S_1,B_1) $  and $\mathfrak{aut}(S_2,B_2)$ into their $\mathfrak{so}(V,g)$-irreducible components.
\end{rema}
\subsection{Tensor Product $S\otimes S$}\hfill
 \vskip 0.1 in
 As usual,  a choice of $B$ on $S$ gives a $C^2$-equivariant isomorphism  $\tau:S\otimes S\rightarrow{ End}(S)$:
$$
\tau(\phi\otimes\psi)(\xi)=B(\phi,\xi)\psi
.$$
Since
$$
B(\tau(\phi\otimes\psi)(\xi),\eta)=(-1)^{\frac{1}{2}n(n-1)}B(\xi,\tau(\psi\otimes\phi)(\eta)),
$$
it follows that
$$
\begin{cases}
\tau(\phi\otimes\psi-\psi\otimes\phi)\in \mathfrak{aut}(S,B)\quad\text{ if }\quad n\equiv 0,1\mod4,\\
\tau(\phi\otimes\psi+\psi\otimes\phi)\in \mathfrak{aut}(S,B)\quad\text{ if }\quad n\equiv 2,3\mod4.
\end{cases}
$$
Identifying ${End}(S)$ with the Clifford algebra $C$ and using the preceding Remark one can reduce symmetric and antisymmetric spinors as $\mathfrak {so}(V,g)$-representations.
\begin{prop}\hfill

 (i) If  $n\equiv 0,1\mod 4$ then $\tau$ induces $\mathfrak{so}(V,g)$-equivariant isomorphisms:
$$
\Lambda^2(S)\cong\bigoplus_{k\equiv 2\text{ or }3\text{ (mod }4)}^{}\hskip-,6cm C^k,\qquad
S^2(S)\cong\bigoplus_{k\equiv 0\text{ or }1\text{ (mod }4)}^{}\hskip-,6cm C^k.
$$

 (ii) If  $n\equiv 2,3\mod4$ then $\tau$ induces $\mathfrak {so}(V,g)$-equivariant isomorphisms:
$$
\Lambda^2(S)\cong\bigoplus_{k\equiv 0\text{ or }1\text{ (mod }4)}^{}\hskip-,6cm C^k,\qquad
S^2(S)\cong\bigoplus_{k\equiv 2\text{ or }3\text{ (mod }4)}^{}\hskip-,6cm C^k.
$$

\end{prop}
 
\section{Cartan's operator $L_2$}
 The operators $L_2$ and $L_{2n}$ to be defined in this section are among the operators $L_k$  that appear in Cartan and Chevalley, where mostly they are used to characterise pure spinors. The interesting properties of $L_2$ to be described herein appear to be new.

Composing $\tau:S\otimes S\rightarrow { End}(S)$ with the projection  $\pi_2:C\rightarrow C^2$
 we can define a $C^2$-equivariant map $L_2:S\times S\rightarrow C^2$ :
 \begin{equation}\label{Ldef}
L_2(\phi,\psi)=\frac{1}{2^{n}}\sum_{i<j}g(e_{i},e_{i})g(e_{j},e_{j})B(\phi,e_je_i\cdot\psi)e_{i}e_{j},
\end{equation}
where $\{e_1,\dots ,e_{2n}\}$ is any orthonormal basis of $V$.

\begin{prop} \label{Lsymmetries}\hfill

\noindent  (i)  If $n\equiv 0$ mod $4$ then $L_2$ is antisymmetric and even, i.e. $L_2 (S_1,S_2) = 0$.

\noindent (ii)  If $n\equiv 1$ mod $4$ then $L_2$ is antisymmetric and odd, i.e. $L_2 (S_1,S_1) = L_2(S_2,S_2) = 0$.

\noindent(iii) If $n\equiv 2$ mod $4$ then $L_2$ is symmetric  and even.

\noindent (iv) If $n\equiv 3$ mod $4$ then $L_2$ is symmetric and odd.
\end{prop}
$L_2$  has an interesting formulation in terms of orbit maps.

\begin{defi} For $\phi\in S$ we define $\phi:V\rightarrow S$ and $\phi^*:S\rightarrow V$ by
\begin{equation}
\phi(v)=v\cdot\phi,\qquad \phi^*(\psi)=\sum_{i}g(e_{i},e_{i})B(\phi,e_i\cdot\psi)e_{i}
\end{equation}
where $\{e_1,\dots ,e_{2n}\}$ is any orthonormal basis of $V$.
\end{defi}
The maps $\phi$ and $\phi^*$ are adjoints for the respective norms, i.e.,
\begin{equation}\label{adjoints}
B(\phi(v),\psi)=g(v,\phi^*(\psi)).
\end{equation}
Thus we have
$$\begin{cases}
g(\phi^*\circ\psi(v_1),v_2)=B(\phi(v_2),\psi(v_1)),\\
 g(\psi^*\circ\phi(v_2),v_1)=B(\psi(v_1),\phi(v_2)),
\end{cases}
$$
hence 
$$
g(\phi^*\circ\psi(v_1),v_2)=(-1)^{\frac{1}{2}n(n-1)}g(\psi^*\circ\phi(v_2),v_1).
$$
Similarly, equation \eqref{adjoints} implies
$$
B(\psi\circ\phi^*(\xi),\eta)=g(\phi^*(\xi),\psi^*(\eta))=B(\phi\circ\psi^*(\eta),\xi),
$$
hence
$$
B(\psi\circ\phi^*(\xi),\eta)=(-1)^{\frac{1}{2}n(n-1)}B(\xi,\phi\circ\psi^*(\eta)).
$$
The next Proposition follows immediately from the above.
\begin{prop} \label{adjsym}
Let $\phi,\psi\in S$ be spinors. Then
$$
(-1)^{\frac{1}{2}n(n-1)}\psi^*\circ\phi-\phi^*\circ\psi,\qquad
(-1)^{\frac{1}{2}n(n-1)}\psi^*\circ\phi+\phi^*\circ\psi
$$
are respectively antisymmetric and symmetric endomorphisms of $(V,g)$, while
$$
(-1)^{\frac{1}{2}n(n-1)}\psi\circ\phi^*-\phi\circ\psi^*,\qquad
(-1)^{\frac{1}{2}n(n-1)}\psi\circ\phi^*+\phi\circ\psi^*
$$
are respectively antisymmetric and symmetric endomorphisms of $(S,B)$.
\end{prop}
Consequently the first expression of Proposition \ref{adjsym} defines a bilinear form on $S$ with values in $\mathfrak{so}(V,g)$, and even/odd depending on $n$ mod $4$.  By Proposition \ref{Lsymmetries},  the bilinear form $L_2$ also takes values in $\mathfrak{so}(V,g)$ and  can be seen to have the same parity properties.  Analogously, the third expression also defines a bilinear form on $S$ with the same parity properties as $L_2$ but taking 
values in $\mathfrak{aut}(S,B)$.  Taken together, they suggest a type of curvature operator.  The exact relationship between them will be described next.

In fact it will be convenient to `renormalise' $L_2$ as follows:
\begin{defi} 
Define $\tilde L_2:S\times S\rightarrow C^2$ by
$$
\tilde L_2=2^{n-1}L_2
$$
\end{defi}

\begin{prop} \label{relations}\hfill

(i) For all $\phi,\psi\in S$,
$$
2\tilde L_2(\phi,\psi)=(-1)^{\frac{1}{2}n(n-1)}\psi^*\circ\phi-\phi^*\circ\psi
$$
in the sense that for all $v\in V$,
\begin{equation}\label{Lid}
2[\tilde L_2(\phi,\psi),v]=(-1)^{\frac{1}{2}n(n-1)}\psi^*\circ\phi(v)-\phi^*\circ\psi(v).
\end{equation}

(ii) For all $\phi,\psi\in S$ and all $v\in V$,
\begin{equation}\label{relation2}
2B(\phi,\psi)v=(-1)^{\frac{1}{2}n(n-1)}\psi^*\circ\phi(v)+\phi^*\circ\psi(v).
\end{equation}
\end{prop}
\begin{demo} Let $\{e_1,\cdots, e_{2n}\}$ be an orthonormal basis of $V$. For any  $e_k$ in this basis,  by equation \eqref{Ldef} we have
$$
[\tilde L_2(\phi,\psi),e_k]=\frac{1}{2}\sum_{i<j}g(e_{i},e_{i})g(e_{j},e_{j})B(\phi,e_je_i\cdot\psi)[e_{i}e_{j},e_k],
$$
and using
$$
[e_ie_j,e_k]=
\begin{cases}
\nonumber
0\qquad\text{if }i\not=k,j\not=k,\\
\nonumber
-2g(e_k,e_k) e_j\quad\text{if }i=k,\\
2g(e_k,e_k) e_i\quad\text{if }j=k,
\end{cases}
$$
this simplifies to
$$
[\tilde L_2(\phi,\psi),e_k]=\sum_{i\not=k}g(e_{i},e_{i})g(e_{k},e_{k})^2B(\phi,e_ke_i\cdot\psi)e_i,
$$
which, since   $g(e_{k},e_{k})^2=1$, reduces to
\begin{equation}\label{eqL1}
[\tilde L_2(\phi,\psi),e_k]=\sum_{i\not=k}g(e_{i},e_{i})B(\phi,e_ke_i\cdot\psi)e_i.
\end{equation}
To calculate the RHS of equation \eqref{Lid} acting on $e_k$ we have 
\begin{align}\label{eqL2}
\nonumber
(-1)^{\frac{1}{2}n(n-1)}\psi^*\circ\phi(e_k)&=(-1)^{\frac{1}{2}n(n-1)}\sum_{i}g(e_{i},e_{i})B(\psi,e_ie_k\cdot\phi)e_{i}\\
\nonumber
&=\sum_{i}g(e_{i},e_{i})B(e_ie_k\cdot\phi,\psi)e_{i}\\
&=\sum_{i}g(e_{i},e_{i})B(\phi,e_ke_i\cdot\psi)e_{i},
\end{align}
and
\begin{equation}\label{eqL3}
\phi^*\circ\psi(e_k)=\sum_{i}g(e_{i},e_{i})B(\phi,e_ie_k\cdot\psi)e_{i}.
\end{equation}
Since $e_k$ is an arbitrary basis element, both parts of the proposition follow from equations \eqref{eqL1}, \eqref{eqL2} and  \eqref{eqL3}.

\end{demo}

In the same way there is a $C^2$-equivariant map $L_{2n}:S\times S\rightarrow C^{2n}$ obtained by composing  $\tau:S\otimes
S\rightarrow {\rm End}(S)=C$ and $\pi_{2n}:C\rightarrow C^{2n}.$  Explicitly,
\begin{equation}\label{defi:L}
L_{2n}(\phi,\psi)=\frac{1}{2^{n}} g(e_{1},e_{1})\dots g(e_{2n},e_{2n})B(\phi,e_{2n}\dots e_1\cdot\psi)e_{1}\dots e_{2n},
\end{equation}
where $\{e_1,\dots ,e_{2n}\}$ is any orthonormal basis of $V$.
The symmetry properties of $L_{2n}$ follow readily from the preceding. 
\begin{prop}\label{ L_{2n}symmetries}\hfill

\noindent(i) If $n\equiv 0$ mod $4$ then $L_{2n}$ is symmetric and even.

\noindent(ii) If $n\equiv 1$ mod $4$ then $L_{2n}$ is antisymmetric and odd.

\noindent(iii) If $n\equiv 2$ mod $4$ then $L_{2n}$ is symmetric  and even.

\noindent(iv) If $n\equiv 3$ mod $4$ then $L_{2n}$ is antisymmetric and odd.
\end{prop}

\subsection{Graded spinor norms}\hfill
\vskip 0.1in
Spinor norms are invariant under $C^2$ but not under $C^1$, as follows from Corollary \ref{spinornorminvprops}.
In order to get something invariant under the action of the Lie algebra $C^1\oplus C^2$ one can use a grading element.

\begin{prop}
Let $B$ be a spinor norm and $\varepsilon\in C$ be a grading element. Define the associated
graded spinor norm $B_{\varepsilon}:S\times S\rightarrow { k}$  by 
$$
B_{\varepsilon}(\phi,\psi)=B(\varepsilon\cdot\phi,\psi)\,\forall \phi,\psi\in S. 
$$
Then
$$
B_{\varepsilon}(v\cdot\phi,\psi)=-B_{\varepsilon}(\phi,v\cdot\psi)\,\,\forall v\in V,\,\,\forall \phi,\psi\in S,
$$
and graded spinor norms are characterised by this property.
\end{prop}

If $\phi,\psi\in S$, $1\le k\le 2n$ and $v_1,\dots ,v_k\in V$ are orthogonal, it is clear that
$$
B_{\varepsilon}(v_1\dots v_k\cdot\phi,\psi)=(-1)^{\frac{1}{2}k(k+1)}B_{\varepsilon}(\phi,v_1\dots v_k\cdot\psi)
$$
and hence we have

\begin{cor}
A graded spinor norm is invariant for the action of the Lie algebra
$C^1\oplus C^2$.  
\end{cor}

\begin{prop}\label{gradedspinornormsymmetries}
 Let $B_{\varepsilon}$ be a graded spinor norm.

\noindent(i) If $n\equiv 0$ mod $4$ then $B_{\varepsilon}$ is symmetric and even.

\noindent(ii) If $n\equiv 1$ mod $4$ then $B_{\varepsilon}$ is antisymmetric and odd.

\noindent(iii) If $n\equiv 2$ mod $4$ then $B_{\varepsilon}$ is antisymmetric  and even.

\noindent(iv) If $n\equiv 3$ mod $4$ then $B_{\varepsilon}$ is symmetric and odd.
\end{prop}

\begin{cor}
 If $\phi,\psi\in S$, $1\le k\le 2n$ and $v_1,\dots ,v_k\in V$ are orthogonal, then
$$
B_{\varepsilon}(\phi,v_1\dots v_k\cdot\psi)=(-1)^{\frac{1}{2}n(n+1)+\frac{1}{2}k(k+1)}B_{\varepsilon}(\psi,v_1\dots v_k\cdot\phi).
$$
\end{cor}

Using $B_{\varepsilon}$ we can define a $C^1\oplus C^2$-equivariant map $\tau_{\varepsilon}:S\times S\rightarrow{\rm End}(S)$ by
$$
\tau_{\varepsilon}(\phi,\psi)(\xi)=B_{\varepsilon}(\phi,\xi)\psi.
$$
One can now repeat the preceding but using the graded versions. 

Using $B_{\varepsilon}$ we can define a $C^1\oplus C^2$-equivariant map $L_{\varepsilon}:S\times S\rightarrow C^1\oplus C^2$ by composing
$\tau_{\varepsilon}:S\otimes S\rightarrow {\rm End}(S)=C$ and $\pi_1\oplus\pi_2:C\rightarrow C^1\oplus C^2.$  (To avoid excessive notation we suppress the subscripts $1,2$.) Explicitly,
$$
L_{\varepsilon}(\phi,\psi)=\frac{1}{2^{n}}\left(\sum_{i}g(e_{i},e_{i})B_{\varepsilon}(\phi,e_i\cdot\psi)e_{i} + \sum_{i<j}g(e_{i},e_{i})g(e_{j},e_{j})B_{\varepsilon}(\phi,e_je_i\cdot\psi)e_{i}e_{j}
\right),
$$
where $\{e_1,\dots ,e_{2n}\}$ is any orthonormal basis of $V$.

\begin{prop}\hfill

\noindent(i) If $n\equiv 0$ mod $4$ then $L_{\varepsilon}$ is antisymmetric, $\pi_2\circ L_{\varepsilon}$ is even and $\pi_1\circ
L_{\varepsilon}$ is odd.

\noindent(ii) If $n\equiv 1$ mod $4$ then $L_{\varepsilon}$ is symmetric, $\pi_2\circ L_{\varepsilon}$ is odd and $\pi_1\circ
L_{\varepsilon}$ is even.

\noindent(iii) If $n\equiv 2$ mod $4$ then $L_{\varepsilon}$ is symmetric, $\pi_2\circ L_{\varepsilon}$ is even and $\pi_1\circ
L_{\varepsilon}$ is odd.

\noindent(iv) If $n\equiv 3$ mod $4$ then $L_{\varepsilon}$ is antisymmetric, $\pi_2\circ L_{\varepsilon}$ is odd and $\pi_1\circ
L_{\varepsilon}$ is even.
\end{prop}

\section{Maximal isotropic subspaces and polarisations }

A maximal isotropic subspace of $V$ is an $n$-dimensional subspace ${\mathcal I}$ of $V$ such that the restriction of $g$ to ${\mathcal I}$
vanishes. By Witt's theorem the group $O(V,g)$ acts transitively on the collection of maximal isotropic subspaces. The stabiliser of ${\mathcal I}$,   $\,S({\mathcal I}),$  is a maximal parabolic subgroup of $O(V,g)$ (see e.g. \cite {Wo}). The natural map from $S({\mathcal I})$ to
$GL({\mathcal I})$ is surjective, giving rise to the exact sequence of groups
$$
1\rightarrow A\rightarrow S({\mathcal I})\rightarrow GL({\mathcal I})\rightarrow 1.
$$
A description of $A$ can be obtained as follows. 

If $a\in A$, then $(a-Id_V)(V)\subseteq {\mathcal I}$, hence there is a unique $\alpha:V\rightarrow V$ such that
\begin{align}
&(i)\quad a=Id_V+\alpha\\
&(ii)\quad \alpha\vert_{\mathcal I}=0,\,\alpha(V)\subseteq {\mathcal I} \text{ and }g(\alpha(v),w)+g(v,\alpha(w))=0.
\end{align}
This shows that $A$ is abelian and $a\mapsto\alpha$ identifies it with a subgroup of the vector space $Hom(V/{\mathcal I},{\mathcal I})$.

 Since ${\mathcal I}$ is maximal isotropic the metric defines a duality pairing $V/{\mathcal I}\otimes {\mathcal I}\rightarrow {k}$, and hence there is an $S({\mathcal I})$-equivariant isomorphism
$$
Hom(V/{\mathcal I},{\mathcal I})\cong {\mathcal I}\otimes {\mathcal I}.
$$

One checks that the composition of maps $A\hookrightarrow Hom(V/{\mathcal I},{\mathcal I})\cong {\mathcal I}\otimes {\mathcal I}$ is an $S({\mathcal I})$-equivariant
isomorphism of $A$ with $\Lambda^2({\mathcal I})$, the space of antisymmetric two tensors on ${\mathcal I}$. Notice that $A$ acts trivially on
$A, \,{\mathcal I},$ and $V/{\mathcal I}$ so that this factors to a $S({\mathcal I})/A\cong GL({\mathcal I})$-equivariant isomorphism 
\begin{equation}
A\cong\Lambda^2({\mathcal I}), 
\end{equation}
as a module for $GL({\mathcal I})$.
Finally, under this
isomorphism the cone of decomposable elements in $\Lambda^2({\mathcal I})$ is the image of the set $T$ of elements $a$ of $A$  with the
property that there exist $v,w\in {\mathcal I}$ such that
$$
a(x)=x+g(x,v)w-g(x,w)v\quad\forall x\in V.
$$
Alternatively, $T$ is  the set of elements of $A$ which are the identity on some $2n-2$ dimensional subspace of $V$
containing ${\mathcal I}$. The subspace determines the group element essentially uniquely and then restriction of the above isomorphism to
$T$ corresponds to the Pl\" ucker embedding of its orthogonal  in $\Lambda^2({\mathcal I}).$

We can realise $S({\mathcal I})$ as a group of affine transformations of an affine space of
which $A$ is the group of translations. Consider the exact sequence of vector spaces:
$$
0\rightarrow {\mathcal I}\rightarrow V \xrightarrow{p} V/{\mathcal I}\rightarrow 0,
$$
and let
\begin{align}
\nonumber
{\mathcal A}&=\{s\in Hom(V/{\mathcal I},V):\,p\circ s=Id_{V/{\mathcal I}}\hbox{ and }\rm{Im}\,s\hbox{ is maximal isotropic}\}\\
\nonumber
&=\{s\in Hom(V/{\mathcal I},V):\, g(s\circ p(v),i)=g(v,i)\,\forall v\in V,i\in {\mathcal I}\text{ and }
\rm{Im}\,s\text{ is maximal isotropic}\}
\end{align}
be the space of isotropic splittings. (One can also identify ${\mathcal A}$ with the space of maximal isotropic complements of ${\mathcal I}$ by $s\in
{\mathcal A}\mapsto Im\,s$). 

This is not a linear subspace of $Hom(V/{\mathcal I},V)$ but it is stable under the natural action
of $S({\mathcal I})$:
$$
f\mapsto g\circ f\circ g^{-1}
$$
and then the group $A$ acts on ${\mathcal A}$ as follows: if $s\in{\mathcal A}$ and if we write $a\in A$ as $a=Id_V+\alpha$, then $a$
maps
$s$ to $s'=s+\alpha$. To check that $s'\in {\mathcal A}$ first note that $p\circ s'=p\circ s=Id_{V/{\mathcal I}}$ and so $s'$ is a
splitting. Further,
\begin{align}
g(s'\circ p(v),s'\circ p(w))&=g(s\circ p(v),\alpha\circ p(w))+g(\alpha\circ p(v),s\circ p(w))\\
&=g(v,\alpha\circ p(w))+g(\alpha\circ p(v),w)\\
&=g(v,\alpha(w))+g(\alpha(v),w)=0,
\end{align}and hence $s'\in {\mathcal A}$. It is clear that $s\mapsto s'$ defines a free group action of $A$ on ${\mathcal A}$. 

To see that $A$ acts transitively, take $s,s'\in{\mathcal A}$. Then $\alpha=s'-s\in Hom(V/{\mathcal I},{\mathcal I})$ and
\begin{align}
0=g(s'\circ p(v),s'\circ p(w))&=g(s\circ p(v),\alpha\circ p(w))+g(\alpha\circ p(v),s\circ  p(w))\\
&=g(v,\alpha\circ p(w))+g(\alpha\circ p(v),w),
\end{align}
and the difference $s'-s$ of the two isotropic splittings is in $A.$

With respect to this affine structure the group $S({\mathcal I})$ acts by affine transformations on ${\mathcal A}$. The tangent space at any
point of ${\mathcal A}$ is canonically isomorphic to $A$ and hence carries a cone structure induced by $T\subset A$ and this is
clearly preserved by the affine action of $S({\mathcal I}).$ From the theory of 3-graded Lie algebras one can see that $S({\mathcal I})$ is exactly the group of affine
transformations of ${\mathcal A}$ which preserve this cone structure.

\begin{rema} Of course for vector spaces $V$ defined over $\mathbb R$ or $\mathbb C$ some of the above is standard (e.g. \cite{Wo}).
\end{rema}
Let ${\mathcal I}$ be a fixed maximal isotropic subspace of $V$. An ${\mathcal I}$-polarisation of $(V,g)$ is a decomposition  $V={\mathcal I}\oplus {\mathcal E}$ such that  ${\mathcal E}$ is a maximal isotropic subspace of $V$, i.e. an element of $\mathcal A$ above.
Choose a  basis $i_1,\dots,i_n$ of ${\mathcal I}$ and let $e_1,\dots,e_n$ be the  basis of ${\mathcal E}$ satisfying
$$
g(e_a,i_b)=\frac{1}{2}\delta_{ab}.
$$
The Clifford
algebra relations
$xy+yx=2g(x,y)Id$ imply that for all $1\le a,b\le n$,
\begin{align}\label{Cliffrel}
\nonumber
&i_ai_b+i_bi_a=0\\
\nonumber
&e_ae_b+e_be_a=0\\
&i_ae_b+e_bi_a=\delta_{ab}.
\end{align}
In terms of this basis the grading operator has the following expression:
\begin{prop}
$$
\varepsilon=(i_1-e_1)(i_1+e_1)\dots (i_n-e_n)(i_n+e_n), \hbox{ in particular }\varepsilon\in C^{2n}.
$$
\end{prop}
For a maximal isotropic subspace ${\mathcal I}$, an ${\mathcal I}$-polarisation $V={\mathcal I}\oplus {\mathcal E}$, and bases $\{i_a\},\{e_a\}$ of ${\mathcal I}$ and ${\mathcal E}$ as above, since $g(i_a\pm e_a,i_b\pm e_b)=\pm \delta_{ab}$ and $g(i_a+ e_a,i_b- e_b)=0$, the set $\{i_a\pm e_a:\,1\le a\le n\}$ is an
orthonormal basis of $V$ as used in \S1. For $1\le a \le n$,  set
$$
E_a=e_a+i_a, E_{\bar a}=e_a-i_a.
$$
Then $\{E_1,E_{\bar 1},\dots,E_n,E_{\bar n}\}$ is  an orthonormal basis with $g(E_a,E_a)=1$ and $g(E_{\bar a},E_{\bar a})=-1$. We order this basis by
$$
1<\bar 1<2<\bar 2<\dots n<\bar n.
$$
The use of $e_i$ in two different ways in a basis, as in \S1 and here in \S5, hopefully does not lead to confusion.

\section{Pure spinors}

\' E. Cartan found a beautiful relationship between the maximal isotropic subspaces of $V$ and a distinguished subset of spinors in $S$. More precisely, he showed that to each maximal isotropic subspace ${\mathcal I}$ there is a unique (up to scalar multiplication) nonzero 
element ${\bf v}_{{\mathcal I}}\in S$ such that
$$
i_a\cdot{\bf v}_{{\mathcal I}}=0\quad\forall 1\le a\le n.
$$
In the language of \'E. Cartan ${\bf v}_{{\mathcal I}}$ is called a pure spinor, and in the language of physics a
vacuum. 

Take an ${\mathcal I}$-polarisation of $(V,g)$, so that $V={\mathcal I}\oplus {\mathcal E}$, and with bases as in \S 4. Then the grading operator, $\varepsilon\in C$, defined there determines the spaces of half-spinors
 $$
\varepsilon\cdot\psi_{\pm}=\pm\psi_{\pm}\quad\forall \psi_{\pm}\in S_{\pm}.
$$

\noindent Repeated use of (\ref{Cliffrel}) shows  that ${\bf v}_{{\mathcal I}}$ is in the half-spinor space $S_+$ associated to $\varepsilon$. 

A basis for $S$ is obtained by applying succesive ``creation operators'' $e_a$ to the ``vacuum'' ${\bf v}_{{\mathcal I}}$, so that $S_+$ is then the
space of ``even particle states'', $S_-$ the space of ``odd particle states'':
\begin{align}
S_+&=<{\bf v}_{{\mathcal I}},e_{i_1}e_{i_2}\dots e_{i_k}\cdot{\bf v}_{{\mathcal I}}:\,1\le i_1<i_2<\dots i_k\le n,\hbox{ k is even}>\cr
S_-&=<e_{i_1}e_{i_2}\dots e_{i_k}\cdot{\bf v}_{{\mathcal I}}:\,1\le i_1<i_2<\dots i_k\le n,\hbox{ k is odd}>.\cr
\end{align}
Notice that $e_{i_1}e_{i_2}\dots e_{i_k}\cdot{\bf v}_{{\mathcal I}}$ is a pure spinor with
$$
{\rm Ann}_V(e_{i_1}e_{i_2}\dots e_{i_k}\cdot{\bf v}_{{\mathcal I}})=<e_{p}, i_{q}: \,p\in\{i_1,\cdots, i_k\},\,q\not\in\{i_1,\cdots, i_k\}>,
$$
a maximal isotropic subspace.
\begin{rema} The basis of pure spinors is an effective computational tool in spinor algebra. Essentially, its use converts computations to combinatorial statements about the parameters for pure spinors.
\end{rema}

An easy example of this is the matrix of the spinor norm $B$. It is helpful to compare this to Proposition \ref{spinornormsymmetries}.
\begin{prop}\label{matrixB}
 If $K=\{k_1,\dots , k_p\}$ and $J=\{j_1,\dots , j_q\}$ are ordered subsets of $\{1,2,\dots
,n\}$  we set $e_K=e_{k_1}\dots e_{k_p}$, $e_J=e_{j_1}\dots e_{j_q}$ and
$e_{\emptyset}=1$. Then
\begin{align}
&(a)\quad B(e_K\cdot{\bf v}_{{\mathcal I}},e_J\cdot{\bf v}_{{\mathcal I}})\not=0\quad 
\Rightarrow 
\quad K\cap J=\emptyset\text{ and }K\cup J=\{1,2,\dots ,n\}\, i.e. \, J=K^c, \\
&(b)\quad B({\bf v}_{{\mathcal I}}, e_1e_2e_3\dots e_n\cdot{\bf v}_{{\mathcal I}})\not= 0.
\end{align}
\end{prop}

For another easy example, take an ${\mathcal I}$-polarisation  $V={\mathcal I}\oplus {\mathcal E}$ and define the following element of the Lie algebra $C^2$:

\begin{equation}\label{3gradopassociatedtopolarisation}
H=\frac{1}{2}\sum_{1}^{n}(e_ai_a-i_ae_a).
\end{equation}

This is independent of the choice of bases $\{i_a\},\{e_a\}$ above. It is an element of $C^2$ since the set
$\{e_ae_b,i_ai_b:1\le a<b\le n\}\cup\{e_ai_b-i_be_a:1\le a,b\le n\}$ is a basis of $C^2.$
\begin{prop}
 For $1\le a \le n$ and $1\le i_1<\dots <i_k\le n$,
$$
(a)\quad [H,e_a]=e_a;\quad (b)\quad [H,i_a]=-i_a;\quad (c)\quad He_{i_1}e_{i_2}\dots
e_{i_k}\cdot{\bf v}_{{\mathcal I}}=(k-\frac{n}{2})e_{i_1}e_{i_2}\dots e_{i_k}\cdot{\bf v}_{{\mathcal I}}.
$$
\end{prop}
$H$ is a useful substitute for  what is called the  number operator in physics, $N=\sum_{a=1}^{a=n}e_ai_a$, which is not in $C^2.$

\section{ Computations of $\tilde L_2$}

\subsection{ The operator $\tilde L_2$}\hfill

Recall the definition of the operator $L_2$:
$$L_2(\phi,\psi)=\frac{1}{2^{n}}\sum_{i<j}g(e_{i},e_{i})g(e_{j},e_{j})B(\phi,e_je_i\cdot\psi)e_{i}e_{j}.$$

Using the basis $\{E_a, E_{\bar b}\}$ from \S 4 we obtain an alternative expression for $\tilde L_2 = 2^{n-1} L_2$ which is more convenient for the computation of the matrix of $\tilde L_2$.
\begin{prop} In terms of the basis of $C^2:$ 
$\{e_ae_b,i_ai_b:1\le a<b\le n\}\cup\{e_ai_b-i_be_a:1\le a,b\le n\}$  we have
\begin{align}
\nonumber
\tilde L_2(\psi_1,\psi_2)&=\sum_{a\not= b}B(e_ae_b\cdot\psi_1,\psi_2)i_ai_b+
\sum_{a\not= b}B(i_ai_b\cdot\psi_1,\psi_2)e_ae_b
+\sum_{a\not= b}B(e_ai_b\cdot\psi_1,\psi_2)(i_ae_b-e_bi_a)\\
\label{Linisotcoords}
&+
\frac{1}{2}\sum_{a}B((e_ai_a-i_ae_a)\cdot\psi_1,\psi_2)(i_ae_a-e_ai_a).
\end{align}
\end{prop}
\begin{demo} 
In a basis $\{E_a, E_{\bar b}\}$
$$
\tilde L_2(\psi_1,\psi_2)=\frac{1}{2}\sum_{i<j}g(E_{i},E_{i})g(E_{j},E_{j})B(\psi_1,E_jE_i\cdot\psi_2)E_{i}E_{j}
$$
and using the ordered basis this sum can be split into two  subsums:
\begin{equation}\label{subsum1}
\frac{1}{2}\sum_{i}g(E_{i},E_{i})g(E_{\bar i},E_{\bar i})B(\psi_1,E_{\bar i}E_i\cdot\psi_2)E_{i}E_{\bar i}
\end{equation}
and 
\begin{align}\label{subsum2}
\nonumber
\frac12\sum_{i<j}
\Big( g(E_{i},E_{i})g(E_{j},E_{ j})B(\psi_1,E_{j}E_i\cdot\psi_2)E_{i}E_{j}\\
\nonumber
+g(E_{i},E_{i})g(E_{\bar j},E_{\bar j})B(\psi_1,E_{\bar j}E_i\cdot\psi_2)E_{i}E_{\bar j}\\
\nonumber
+g(E_{\bar i},E_{\bar i})g(E_{j},E_{j})B(\psi_1,E_{j}E_{\bar i}\cdot\psi_2)E_{\bar i}E_{j}\\
+g(E_{\bar i},E_{\bar i})g(E_{\bar j},E_{\bar j})B(\psi_1,E_{\bar j}E_{\bar i}\cdot\psi_2)E_{\bar i}E_{\bar j} \Bigr).
\end{align}
Since
$$
E_aE_{\bar a}=-e_ai_a+i_ae_a,\quad E_{\bar a}E_a=-i_ae_a+e_ai_a
$$
the sum \eqref{subsum1} reduces to
\begin{equation}\label{}
-\frac{1}{2}\sum_{a}B(\psi_1,(e_ai_a-i_ae_a)\cdot\psi_2)(i_ae_a-e_ai_a)=\frac{1}{2}\sum_{a}B((e_ai_a-i_ae_a)\cdot\psi_1,\psi_2)(i_ae_a-e_ai_a)
\end{equation}
and this is the last term in \eqref{Linisotcoords}.

To simplify \eqref{subsum2} we first observe that
\begin{align}
\nonumber
E_aE_{ b}= e_{ab}+e_ai_b+i_ae_b+i_ai_b,\quad E_aE_{\bar b}=e_ae_b-e_ai_b+i_ae_b-i_ai_b,\\
E_{\bar a}E_{ b}= e_{ab}+e_a i_b-i_ae_b-i_ai_b,\quad E_{\bar a}E_{\bar b}=e_ae_b-e_ai_b-i_ae_b+i_ai_b,
\end{align}
Hence for fixed $a<b$, the coefficient of $e_{ab}$ in \eqref{subsum2} is
$$
\frac{1}{2}
\left(
B(\psi_1,E_{b}E_a\cdot\psi_2)-B(\psi_1,E_{\bar b}E_a\cdot\psi_2)-B(\psi_1,E_{b}E_{\bar a}\cdot\psi_2)+B(\psi_1,E_{\bar b}E_{\bar a}\cdot\psi_2)
\right)
$$
which can be written
$$
\frac{1}{2}B(\psi_1,(E_{b}-E_{\bar b})(E_a-E_{\bar a})\cdot\psi_2),
$$
that is
$$
2B(\psi_1,i_bi_a\cdot\psi_2)=2B(i_ai_b\cdot\psi_1,\psi_2).
$$
Summing over all $a<b$ we get
$$
2\sum_{1\le a<b\le n}B(i_ai_b\cdot\psi_1,\psi_2)e_{ab}
$$
which is the second term in \eqref{Linisotcoords}. Similarly, looking at the coefficients of $i_ai_b$ and $i_ae_b-e_bi_a$ in \eqref{subsum2} we get the first and third terms of \eqref{Linisotcoords}. 
\end{demo}

\subsection{ The matrix of $\tilde L_2$}\hfill

Next we use the basis of pure spinors to simplify the expression. Let ${\bf v}_{{\mathcal I}}$ be a pure spinor defined by the maximal isotropic subspace ${\mathcal I}$. Recall from \S 5 the basis of pure spinors and take $\psi_1 = e_I\cdot{\bf v}_{{\mathcal I}}$ and $\psi_2 = e_J\cdot{\bf v}_{{\mathcal I}}$. 
Most terms in the formula for $\tilde L_2(e_I\cdot{\bf v}_{{\mathcal I}},e_J\cdot{\bf v}_{{\mathcal I}})$ vanish by Proposition \ref{matrixB}:
\begin{align}
\nonumber
&(a)\quad B(e_ae_be_I\cdot{\bf v}_{{\mathcal I}},e_J\cdot{\bf v}_{{\mathcal I}})\not=0\quad \Leftrightarrow \quad I\cap J=\emptyset\text{ and } I^c\cap J^c=\{a,b\}.\\
\nonumber
&(b)\quad B(i_ai_be_I\cdot{\bf v}_{{\mathcal I}},e_J\cdot{\bf v}_{{\mathcal I}})\not=0\quad  \Leftrightarrow \quad I\cap J=\{a,b\}\text{ and } I^c\cap J^c=\emptyset.\\
\nonumber
&(c)\quad B(e_ai_be_I\cdot{\bf v}_{{\mathcal I}},e_J\cdot{\bf v}_{{\mathcal I}})\not=0\text{ and }a\not= b\quad \Leftrightarrow \quad I\cap J=\{b\}\hbox{ and } I^c\cap
J^c=\{a\}.\\ 
\nonumber
&(d)\quad B((e_ai_a -i_ae_a)e_I\cdot{\bf v}_{{\mathcal I}},e_J\cdot{\bf v}_{{\mathcal I}})\not=0\quad  \Leftrightarrow \quad I\cap J=\emptyset \text{ and } I^c\cap
J^c=\emptyset. 
\end{align}
\begin{rema} Note that $ (I\cap J)\cup  (I^c\cap J^c)=(I\Delta J)^c$ where $I\Delta J$ denotes the symmetric difference of the sets $I$ and $J$.
\end{rema}
From this we can calculate the matrix of $\tilde L_2(e_I\cdot{\bf v}_{{\mathcal I}},e_J\cdot{\bf v}_{{\mathcal I}})$ in the basis of particle states $\{ e_K\cdot{\bf v}_{{\mathcal I}}\}$.  
 
\vskip,2cm
\begin{prop}\label{Lnonzero}
$\tilde L_2(e_I\cdot{\bf v}_{{\mathcal I}},e_J\cdot{\bf v}_{{\mathcal I}})\not\equiv0$ iff $I$ and $J$ satisfy one of (a), (b), (c), (d) above. In those cases in terms of the basis  $\{ e_K\cdot{\bf v}_{{\mathcal I}}\}$ we have
\begin{align}
\nonumber
(a) \text{ If  } I\cap J&=\emptyset\text{ and } I^c\cap J^c=\{a,b\} \text{ then} \\
\nonumber
&\tilde L_2(e_I\cdot{\bf v}_{{\mathcal I}},e_J\cdot{\bf v}_{{\mathcal I}})e_K\cdot{\bf v}_{{\mathcal I}}=2B(e_ae_be_I\cdot{\bf v}_{{\mathcal I}},e_J\cdot{\bf v}_{{\mathcal I}})i_ai_b e_K\cdot{\bf v}_{{\mathcal I}}.\\
\nonumber
(b) \text{ If  } I\cap J&=\{a,b\}\text{ and } I^c\cap J^c=\emptyset \text{ then} \\
\nonumber
&\tilde L_2(e_I\cdot{\bf v}_{{\mathcal I}},e_J\cdot{\bf v}_{{\mathcal I}})e_K\cdot{\bf v}_{{\mathcal I}}=2B(i_ai_be_I\cdot{\bf v}_{{\mathcal I}},e_J\cdot{\bf v}_{{\mathcal I}})e_ae_b e_K\cdot{\bf v}_{{\mathcal I}}.\\
\nonumber
(c) \text{ If  }  I\cap J&=\{b\}\text{ and } I^c\cap J^c=\{a\} \text{ then}  \\
\nonumber
&\tilde L_2(e_I\cdot{\bf v}_{{\mathcal I}},e_J\cdot{\bf v}_{{\mathcal I}})e_K\cdot{\bf v}_{{\mathcal I}}=2B(e_ai_be_I\cdot{\bf v}_{{\mathcal I}},e_J\cdot{\bf v}_{{\mathcal I}})i_ae_be_K\cdot{\bf v}_{{\mathcal I}}.\\ 
\nonumber
(d) \text{ If  } I\cap J&=\emptyset \text{ and } I^c\cap J^c=\emptyset \text{ then}  \\
\nonumber
&\tilde L_2(e_I\cdot{\bf v}_{{\mathcal I}},e_J\cdot{\bf v}_{{\mathcal I}})e_K\cdot{\bf v}_{{\mathcal I}}=\frac12B(e_I\cdot{\bf v}_{{\mathcal I}},e_J\cdot{\bf v}_{{\mathcal I}})(n-2\vert I\cap K\vert-2\vert I^c\cap K^c\vert)e_K\cdot{\bf v}_{{\mathcal I}}. 
\end{align}
\end{prop}

\begin{rema} Given two pure spinors $\psi, \psi'$ such that $\tilde L_2(\psi,\psi')\not\not\equiv0$, the intersection properties of the associated maximal isotropic subspaces ${\rm Ann}_V(\psi),{\rm Ann}_V(\psi')$ determine $\tilde L_2(\psi,\psi')$  up to a constant:
\begin{itemize}
\item
if $B(\psi,\psi')=0$ then dim$({\rm Ann}_V(\psi)\cap{\rm Ann}_V(\psi'))=2$ and $\tilde L_2(\psi,\psi')$ is proportional to $Q(\omega)$ for any nonzero  $\omega\in\Lambda^2({\rm Ann}_V(\psi)\cap{\rm Ann}_V(\psi'))$ (see \cite{EC}).
\item
if $B(\psi,\psi')\not=0$ then dim$({\rm Ann}_V(\psi)\cap{\rm Ann}_V(\psi'))=0$ and  $\tilde L_2(\psi,\psi')$ is proportional to the operator $H$ associated to the polarisation $V={\rm Ann}_V(\psi)\oplus{\rm Ann}_V(\psi')$ (see \eqref{3gradopassociatedtopolarisation}).
\end{itemize}
This is a weaker but  `geometric' version of Proposition \ref{Lnonzero}.
For example, if  $\psi=e_I\cdot{\bf v}_{{\mathcal I}}$ and $\psi'=e_J\cdot{\bf v}_{{\mathcal I}}$ are as in  Proposition \ref{Lnonzero} (a), we have 
$$
B(\psi,\psi')=0,\quad {\rm Ann}_V(\psi)\cap{\rm Ann}_V(\psi')=<i_a,i_b>
$$
and this result implies that $\tilde L_2(\psi,\psi')$ is proportional to $Q(i_a\wedge i_b)=i_ai_b$, whereas more importantly Proposition \ref{Lnonzero} also gives the constant of proportionality.
\end{rema}

There is still some simplification possible in the parameters $I,J,K$. Looking at the above more closely we see that $\tilde L_2(e_I\cdot{\bf v}_{{\mathcal I}},e_J\cdot{\bf v}_{{\mathcal I}})e_K\cdot{\bf v}_{{\mathcal I}}$ is `symmetric' in $I,J,K$ in the following sense.

\begin{cor}\hfill

\noindent (i) If either $I\cap J\cap K\not=\emptyset$ or $I^c\cap J^c\cap K^c\not=\emptyset$ then
$$\tilde L_2(e_I\cdot{\bf v}_{{\mathcal I}},e_J\cdot{\bf v}_{{\mathcal I}})e_K\cdot{\bf v}_{{\mathcal I}}=0.$$

\noindent (ii) If $I\cap J\cap K=\emptyset$ and $I^c\cap J^c\cap K^c=\emptyset$ then 
$$\tilde L_2(e_I\cdot{\bf v}_{{\mathcal I}},e_J\cdot{\bf v}_{{\mathcal I}})e_K\cdot{\bf v}_{{\mathcal I}}\textrm{ is proportional to } e_{(I\cap J)\cup(J\cap K)\cup(K\cap I)}\cdot{\bf v}_{{\mathcal I}}.$$
\end{cor}
Now three subsets $I,J,K$ of $\{1,2,\dots, n\}$ satisfying the conditions 
\begin{equation}\label{subsetcondition}
I\cap J\cap K=I^c\cap J^c\cap K^c=\emptyset
\end{equation}
define a partition of
$\{1,2,\dots, n\}$ into six disjoint subsets:
$$
\{1,2,\dots, n\}=(I\cap J)\cup  (J\cap K)\cup (K\cap I )\cup (I^c\cap J^c)\cup  (J^c\cap K^c)\cup (K^c\cap I^c)
$$
and in terms of this  partition
\begin{align}
\nonumber
I&=(I\cap J)\cup  (K\cap I)\cup (J^c\cap K^c) \\
\nonumber
J&=(J\cap K)\cup  (I\cap J)\cup (K^c\cap I^c) \\
K&=(K\cap I)\cup  (J\cap K)\cup (I^c\cap J^c). 
\end{align}
The simplest example of three subsets satisfying the conditions  \eqref{subsetcondition} is given by three pairwise disjoint subsets of  $\{1,2,\dots, n\}$ whose union is $\{1,2,\dots, n\}$. In fact this is the general case.
\begin{prop} 
Let  $I,J,K$ be three oriented subsets of  $\{1,2,\dots, n\}$. Suppose that \hfill 

\noindent $I\cap J\cap K=\emptyset$ and $I^c\cap J^c\cap K^c=\emptyset$.  Then there is a polarisation 
$V={\mathcal I}'\oplus{\mathcal E}'$ and oriented subsets $I',J',K'$  of  $\{1,2,\dots, n\}$ such that 
\vskip0,2cm
(i) $e_I\cdot{\bf v}_{{\mathcal I}}=e'_{I'}\cdot{\bf v}_{{\mathcal I}'}$, \quad$e_J\cdot{\bf v}_{{\mathcal I}}=e'_{J'}\cdot{\bf v}_{{\mathcal I}'}$\quad 
and \quad$e_K\cdot{\bf v}_{{\mathcal I}}=e'_{K'}\cdot{\bf v}_{{\mathcal I}'}$.
\vskip0,2cm
(ii)  $I'\cap J'=\emptyset$, $K'=I'^c\cap J'^c$, and $I'\cup J'\cup K'=\{1,2,\dots, n\}$. 

\vskip0,2cm
(iii) If $I,J,K$ are of the same parity then  $I',J',K'$ are of the same parity, i.e $S_+$ or $S_-$.
\end{prop}
\begin{demo} Set
\begin{align}
\nonumber
{\mathcal I}'&={\rm Vect} <e_a,i_b:\, a\in (I^c\cap J^c)\cup  (J^c\cap K^c)\cup (K^c\cap I^c),b\in
 (I\cap J)\cup  (J\cap K)\cup (K\cap I )>,\\
\nonumber
{\mathcal E}'&={\rm Vect} <e_a,i_b:\,
a\in(I\cap J)\cup  (J\cap K)\cup (K\cap I ),\,
b\in (I^c\cap J^c)\cup  (J^c\cap K^c)\cup (K^c\cap I^c)>.
\end{align}
Then it is clear that $V={\mathcal I}'\oplus{\mathcal E}'$ is  a polarisation, and that 
$$
{\bf v}_{{\mathcal I}'} :=e_{(I\cap J)\cup  (J\cap K)\cup (K\cap I )}\cdot{\bf v}_{{\mathcal I}}
$$
is a  pure spinor defined by ${\mathcal I}'$. It is equally clear that if 
\begin{align}
\nonumber
I'&=(J^c\cap K^c)\cup  (J\cap K), \\
\nonumber
J'&=(K^c\cap I^c)\cup  (K\cap I),\\
\nonumber
K'&=(I^c\cap J^c)\cup  (I\cap J), 
\end{align}
then $I'\cap J'=J'\cap K'=K'\cap I'=\emptyset$ and $I'\cup J'\cup K'=\{1,2,\dots, n\}$. 

Define
$$
e'_a=\
\begin{cases}
\nonumber
e_a\quad \text{ if }\,a\in (I^c\cap J^c)\cup  (J^c\cap K^c)\cup (K^c\cap I^c),\\
\nonumber
i_a\quad\text{ if }\,a \in (I\cap J)\cup  (J\cap K)\cup (K\cap I )
\end{cases}
$$
and
$$
i'_a=\
\begin{cases}
\nonumber
e_a\quad \text{ if }\,a \in (I\cap J)\cup  (J\cap K)\cup (K\cap I ),\\
\nonumber
i_a\quad\text{ if }\,a\in (I^c\cap J^c)\cup  (J^c\cap K^c)\cup (K^c\cap I^c).
\end{cases}
$$
Then $\{e'_a,i'_a:\,1\le a\le n\}$ satisfy the Clifford relations 
\begin{align}
\nonumber
&i'_ai'_b+i'_bi'_a=0\\
\nonumber
&e'_ae'_b+e'_be'_a=0\\
&i'_ae'_b+e'_bi'_a=\delta_{ab}
\end{align}
and up to signs,
$$
e'_{I'}\cdot{\bf v}_{{\mathcal I}'}=e_{J^c\cap K^c}i_{J\cap K}e_{(I\cap J)\cup  (J\cap K)\cup (K\cap I )}\cdot{\bf v}_{{\mathcal I}}=e_{J^c\cap K^c}e_{(I\cap J)\cup (K\cap I )}\cdot{\bf v}_{{\mathcal I}}.
$$
Since $(J^c\cap K^c)\cup{(I\cap J)\cup (K\cap I )}=I$ this means (up to signs)
$$
e'_{I'}\cdot{\bf v}_{{\mathcal I}'}=e_I\cdot{\bf v}_{{\mathcal I}}.
$$
Similarly $e'_{J'}\cdot{\bf v}_{{\mathcal I}'}=e_J\cdot{\bf v}_{{\mathcal I}}$ and $e'_{K'}\cdot{\bf v}_{{\mathcal I}'}=e_K\cdot{\bf v}_{{\mathcal I}}$. This proves (i).  Parts (ii) and (iii) follow immediately.
\end{demo}
\begin{prop}\label{1explicitformforL}\, Let  $I,J,K$ be three oriented subsets of  $\{1,2,\dots, n\}$ that are pairwise disjoint, 
 $I\cap J=J\cap K=K\cap I=\emptyset$, and $I^c\cap J^c\cap K^c=\emptyset$.  Then
\vskip0,2cm
(i) $\vert I\vert +\vert J\vert +\vert K\vert =n$.
\vskip0,2cm
(ii) If $\vert K\vert $ is not equal to $0$ or $2$, then
$\tilde L_2(e_I\cdot{\bf v}_{{\mathcal I}},e_J\cdot{\bf v}_{{\mathcal I}})e_K\cdot {\bf v}_{{\mathcal I}}=0$.
\vskip0,2cm
(iii) If $\vert K\vert =0$ then
$$
\tilde L_2(e_I\cdot{\bf v}_{{\mathcal I}},e_J\cdot{\bf v}_{{\mathcal I}})e_K\cdot {\bf v}_{{\mathcal I}}=
B(e_I\cdot{\bf v}_{{\mathcal I}},e_J\cdot{\bf v}_{{\mathcal I}})(\vert I\vert -\frac n2){\bf v}_{{\mathcal I}}.
$$
\vskip0,2cm
(iv) If $\vert K\vert =2$ and $K=\{\overrightarrow{ba}\}$ then
$$
\tilde L_2(e_I\cdot{\bf v}_{{\mathcal I}},e_J\cdot{\bf v}_{{\mathcal I}})e_K\cdot{\bf v}_{{\mathcal I}}=
2B(e_{ab}e_I\cdot{\bf v}_{{\mathcal I}},e_J\cdot{\bf v}_{{\mathcal I}}){\bf v}_{{\mathcal I}}.
$$
\end{prop}
\begin{demo} 
Part (i) is clear since as we observed before $I,J,K$ define a partition of $\{1,\cdots, n\}$. Parts (ii), (iii) and (iv)  follow from Proposition \ref{Lnonzero} since $I\cap J=\emptyset$ and $K=I^c\cap J^c$.
\end{demo}

In \S 8 we will need an expression for $L_{2n}$. Recall from (\ref{defi:L})

$$L_{2n}(\psi_1,\psi_2)=\frac{1}{2^{n}} g(e_{1},e_{1})\dots g(e_{2n},e_{2n})B(\psi_1,e_{2n}\dots e_1\cdot\psi_2)e_{1}\dots e_{2n}.$$
A simplification of this in terms of pure spinors is rather straightforward.   
\begin{prop}
$$
L_{2n}(\psi_1,\psi_2)=\frac{1}{2^n}B(\psi_1,\varepsilon\cdot\psi_2)\varepsilon.
$$
\end{prop}

\subsection{R\'esum\'e}\hfill

Let
\begin{itemize}
\item $(V,g)$ be a $2n$-dimensional vector space with a hyperbolic metric $g$;
\item $S$ be  a space of spinors  (i.e., we identify $C(V,g)$ with ${\rm End}(S)$ for some
$2^n$-dimensional vector space $S$;
\item $B:S\times S\rightarrow k$ be a spinor norm (Cartan form).
\item
$V={\mathcal I}\oplus {\mathcal E}$ be a polarisation of $V$;
\item
${\bf v}_{{\mathcal I}}$ be a pure spinor associated to ${\mathcal I}$  (i.e., $v\cdot{\bf v}_{{\mathcal I}}=0$ for all $v\in {\mathcal I})$;
\item $i_1,\dots ,i_n$ and $e_1,\dots ,e_n$ be  bases of ${\mathcal I}$ and ${\mathcal E}$ respectively such that
\begin{align}\label{Cliffrel1}
\nonumber
&i_ai_b+i_bi_a=0\\
\nonumber
&e_ae_b+e_be_a=0\\
&i_ae_b+e_bi_a=\delta_{ab};
\end{align}
\item $\tilde L_2:S\times S\rightarrow C^2(V,g)$ be the normalised projection operator:
\begin{align}
\nonumber
\tilde L_2(\psi_1,\psi_2)&=\sum_{a\not= b}B(e_ae_b\psi_1,\psi_2)i_ai_b+
\sum_{a\not= b}B(i_ai_b\psi_1,\psi_2)e_ae_b
+\sum_{a\not= b}B(e_ai_b\psi_1,\psi_2)(i_ae_b-e_bi_a)\\
&+
\frac{1}{2}\sum_{a}B((e_ai_a-i_ae_a)\cdot\psi_1,\psi_2)(i_ae_a-e_ai_a).
\end{align}
\end{itemize}
Then for all oriented subsets $I,J,K$ of $\{1,\dots ,n\}$ we have shown that:
\begin{itemize}
\item  $\tilde L_2(e_I\cdot{\bf v}_{{\mathcal I}},e_J\cdot{\bf v}_{{\mathcal I}})e_K\cdot{\bf v}_{{\mathcal I}}=0$ unless $I\cap J\cap K=\emptyset$ and $I^c\cap J^c\cap K^c=\emptyset$.
\item If $I\cap J\cap K=\emptyset$ and $I^c\cap J^c\cap K^c=\emptyset$ there is a polarisation
$V={\mathcal I}'\oplus {\mathcal E}'$ and oriented subsets $I',J',K'$ of $\{1,\cdots ,n\}$ such that
\begin{align}
\nonumber
&(i)\quad e_I\cdot{\bf v}_{{\mathcal I}}=e_{I'}\cdot{\bf v}_{{\mathcal I}'},\quad
e_J\cdot{\bf v}_{{\mathcal I}}=e_{J'}\cdot{\bf v}_{{\mathcal I}'},\quad
e_K\cdot{\bf v}_{{\mathcal I}}=e_{K'}\cdot{\bf v}_{{\mathcal I}'}.\\
&(ii)\quad I'\cap J'=\emptyset\text{ and }K'=I'^c\cap J'^c.
\end{align}
\begin{prop} \label{explicitformforL}\item
Let  $I,J,K$ be three oriented subsets of  $\{1,2,\dots, n\}$ satisfying $I\cap J=\emptyset$ and $K=I^c\cap J^c$.  Then
\vskip0,2cm
\noindent(i) $\vert I\vert +\vert J\vert +\vert K\vert =n$.
\vskip0,2cm
\noindent(ii) If $\vert K\vert $ is not equal to $0$ or $2$, then
$\tilde L_2(e_I\cdot{\bf v}_{{\mathcal I}},e_J\cdot{\bf v}_{{\mathcal I}})e_K\cdot {\bf v}_{{\mathcal I}}=0$.
\vskip0,2cm
\noindent(iii) If $\vert K\vert =0$ then
$$
\tilde L_2(e_I\cdot{\bf v}_{{\mathcal I}},e_J\cdot{\bf v}_{{\mathcal I}})e_K\cdot {\bf v}_{{\mathcal I}}=
B(e_I\cdot{\bf v}_{{\mathcal I}},e_J\cdot{\bf v}_{{\mathcal I}})(\vert I\vert -\frac n2){\bf v}_{{\mathcal I}}.
$$
\vskip0,2cm
\noindent(iv) If $\vert K\vert =2$ and $K=\{\overrightarrow{ba}\}$ then
$$
\tilde L_2(e_I\cdot{\bf v}_{{\mathcal I}},e_J\cdot{\bf v}_{{\mathcal I}})e_K\cdot{\bf v}_{{\mathcal I}}=
2B(e_{ab}e_I\cdot{\bf v}_{{\mathcal I}},e_J\cdot{\bf v}_{{\mathcal I}}){\bf v}_{{\mathcal I}}.
$$
\end{prop}
\item In all cases $\tilde L_2(e_I\cdot{\bf v}_{{\mathcal I}},e_J\cdot{\bf v}_{{\mathcal I}})e_K\cdot{\bf v}_{{\mathcal I}}$ is proportional to
$e_{(I\cap J)(J\cap K)(K\cap I)}\cdot{\bf v}_{{\mathcal I}}$.
\end{itemize}

\section{ Potential Lie algebra structures}
As Cartan's operator $\tilde L_2$ (or $\tilde L_2 + L_{2n}$) maps from $S_i\times S_i$ (or $S\times S$) to $C^2$ (or $C^2\oplus C^{2n}$) it provides a natural candidate for a type of \lq\lq curvature\rq\rq operator on $S_i$ (or $S$). Cartan's calculation of curvature operators for symmetric spaces then motivates possible Lie triple system structures.

If $n\equiv 0$ mod $4$, we can now define ($i=1,2$) a unique antisymmetric map

$$(C^2\oplus S_i)\times(C^2\oplus S_i)\rightarrow
C^2\oplus S_i$$ such that
\begin{align}
[A,B]&=AB-BA\quad\hbox{ if }A,B\in C^2;\\
[A,\psi]&=A\cdot\psi\quad\hbox{ if }A\in C^2,\,\psi\in S_{i};\\
[\phi,\psi]&=\tilde L_2(\phi,\psi)\quad\hbox{ if }\phi,\psi\in S_{i}.
\end{align}
Similarly, If $n\equiv 1$ mod $4$, we can define  a unique antisymmetric map 
$$(C^2\oplus C^{2n}\oplus S)\times(C^2\oplus
C^{2n}\oplus S)\rightarrow C^2\oplus C^{2n}\oplus S$$ such that
\begin{align}
[A,B]&=AB-BA\quad\hbox{ if }A,B\in C^2\oplus C^{2n}\\
[A,\psi]&=A\cdot\psi\quad\hbox{ if }A\in C^2\oplus C^{2n},\,\psi\in S\\
[\phi,\psi]&=\tilde L_2(\phi,\psi)+L_{2n}(\phi,\psi)\quad\hbox{ if }\phi,\psi\in S.
\end{align}
The question is: do these brackets define Lie algebra structures on
$C^2\oplus S_i\quad  (n=0\hbox{ mod }4)$ and 
$C^2\oplus C^{2n}\oplus S \quad (n=1\hbox{ mod }4)$ respectively? Since $C^2,C^2\oplus C^{2n}$ are Lie algebras , since $S_i,S$
are representations and  since $\tilde L_2:S_i\times S_i\rightarrow C^2$ and $\tilde L_2+L_{2n}:S\times S\rightarrow C^2\oplus C^{2n}$ are 
equivariant maps, this will be the case if and only if the following Jacobi identities are satisfied:
$$
\tilde L_2(\psi_1,\psi_2)\cdot\psi_3+\tilde L_2(\psi_2,\psi_3)\cdot\psi_1+\tilde L_2(\psi_3,\psi_1)\cdot\psi_2=0\quad (n=0\hbox{ mod }4)
$$ 
and respectively if for some $a,b\in k^*$ 
$$
(a\tilde L_2+bL_{2n})(\psi_1,\psi_2)\cdot\psi_3+(a\tilde L_2+bL_{2n})(\psi_2,\psi_3)\cdot\psi_1+(a\tilde L_2+bL_{2n})(\psi_3,\psi_1)\cdot\psi_2=0\quad (n=1\hbox{ mod }4).
$$
Note that if $n=8$, $C^2\oplus S_i$ is of dimension $248$ and that
if $n=5$, $C^2\oplus C^{2n}\oplus S$ is of dimension $78.$

\section{Spinor constructions of exceptional Lie algebras}
\subsection{Construction of split $e_8$}\hfill

Let $(V,g)$ be a sixteen-dimensional vector space with a nondegenerate hyperbolic symmetric bilinear form $g$. Choose a 256-dimensional space of spinors $S=S_1\oplus S_2$ and  an isomorphism  $C(V,g)\cong End(S)$. Since $n=8=0$ mod 4, $B:S\times S\rightarrow k$ is even symmetric and $\tilde L_2:S\times S\rightarrow C^2(V,g)$  is even antisymmetric  (Proposition \ref{spinornormsymmetries} and Proposition \ref{Lsymmetries}).

On  the 248-dimensional vector space
$$
E=C^2(V,g)\oplus S_1
$$
following the procedure above we define $[\phantom x ,\phantom y ]:E\otimes E\rightarrow E$ to be the unique antisymmetric bilinear map such that:
\begin{align}
[A,B]&=AB-BA\quad\hbox{ if }A,B\in C^2;\\
[A,\psi]&=A\cdot\psi\quad\hbox{ if }A\in C^2,\,\psi\in S_1;\\
[\phi,\psi]&=\tilde L_2(\phi,\psi)\quad\hbox{ if }\phi,\psi\in S_1.
\end{align}
where $\tilde L_2$ is given by \eqref{Linisotcoords}.  The bracket $[\phantom x ,\phantom y ]$ defines a Lie algebra structure on $E$ iff 
\begin{equation}\label{Jac1fore_8}
\tilde L_2(\psi_1,\psi_2)\cdot\psi_3+\tilde L_2(\psi_2,\psi_3)\cdot\psi_1+\tilde L_2(\psi_3,\psi_1)\cdot\psi_2=0
\quad\forall \psi_1,\psi_2,\psi_3\in S_1.
\end{equation}
Choose a polarisation $V={\mathcal I}\oplus {\mathcal E}$  of $V$ such that a  pure spinor ${\bf v}_{{\mathcal I}}$
corresponding to ${\mathcal I}$ is in $S_1$ and bases $i_1,\cdots ,i_n$ and $e_1,\cdots ,e_n$ of respectively ${\mathcal I}$ and ${\mathcal E}$ such that
\begin{align}\label{Cliffrel2}
\nonumber
&i_ai_b+i_bi_a=0\\
\nonumber
&e_ae_b+e_be_a=0\\
&i_ae_b+e_bi_a=\delta_{ab}.
\end{align}
Then 
$$
\{e_I\cdot{\bf v}_{{\mathcal I}}:\quad\vert I\vert\text{ is even or $\emptyset$}\}
$$
is a basis of $S_1$ and to prove the Jacobi identity \eqref{Jac1fore_8} it is sufficient to prove that
\begin{equation}\label{Jac2fore_8}
\tilde L_2(e_I\cdot{\bf v}_{{\mathcal I}},e_J\cdot{\bf v}_{{\mathcal I}})\cdot e_K\cdot{\bf v}_{{\mathcal I}}+
\tilde L_2(e_J\cdot{\bf v}_{{\mathcal I}},e_K\cdot{\bf v}_{{\mathcal I}})\cdot e_I\cdot{\bf v}_{{\mathcal I}}+
\tilde L_2(e_K\cdot{\bf v}_{{\mathcal I}},e_I\cdot{\bf v}_{{\mathcal I}})\cdot e_J\cdot{\bf v}_{{\mathcal I}}
=0
\end{equation}
for all even subsets $I,J,K$ of $\{1,2,3,4,5,6,7,8\}$. 

By the r\'esum\'e above,  all terms in this equation vanish unless $I\cap J\cap K=\emptyset$ and $I^c\cap J^c\cap K^c=\emptyset$. So in fact we need only prove \eqref{Jac2fore_8} for even subsets of $\{1,2,3,4,5,6,7,8\}$ satisfying these two conditions. But then, again by the r\'esum\'e, by changing the polarisation if necessary, we can always assume that
$$
I\cap J=\emptyset,\quad K=I^c\cap J^c
$$
and then all terms in \eqref{Jac2fore_8} vanish unless one of the sets $I,J,K$ has $0$ or $2$ elements. 
So in the end,  to prove that the bracket $[\phantom x ,\phantom y ]$ defines a Lie algebra structure on $E$ it  remains to prove only that
\begin{equation}\label{Jac22fore_8}
\tilde L_2(e_I\cdot{\bf v}_{{\mathcal I}},e_J\cdot{\bf v}_{{\mathcal I}})\cdot e_K\cdot{\bf v}_{{\mathcal I}}+
\tilde L_2(e_J\cdot{\bf v}_{{\mathcal I}},e_K\cdot{\bf v}_{{\mathcal I}})\cdot e_I\cdot{\bf v}_{{\mathcal I}}+
\tilde L_2(e_K\cdot{\bf v}_{{\mathcal I}},e_I\cdot{\bf v}_{{\mathcal I}})\cdot e_J\cdot{\bf v}_{{\mathcal I}}
=0
\end{equation}
for those  oriented subsets $I,J,K$ of $\{1,2,3,4,5,6,7,8\}$ such that
\begin{itemize}
\item  $\vert I\vert, \vert J\vert$ and $\vert K\vert$ are even;
\item $\vert I\vert+\vert J\vert+\vert K\vert=8$;
\item $I\cap J= J\cap K= K\cap I=\emptyset$;
\item One of $\vert I\vert, \vert J\vert$ or $\vert K\vert$ is equal to $0$ or $2$.
\end{itemize}
Up to permutations of $I,J$ and $K$ the only possibilities are 
\begin{itemize}
\item[(i)]  $\vert I\vert=0$, $\vert J\vert=0$ and $\vert K\vert=8$;
\item[(ii)] $\vert I\vert=0$, $\vert J\vert=2$ and $\vert K\vert=6$;
\item[(iii)]  $\vert I\vert=0$, $\vert J\vert=4$ and $\vert K\vert=4$;
\item[(iv)]  $\vert I\vert=2$, $\vert J\vert=2$ and $\vert K\vert=4$.
\end{itemize}
Since $\vert K\vert\not=0$ and $\vert K\vert\not=2$ in all four cases we have
$$
\tilde L_2(e_I\cdot{\bf v}_{{\mathcal I}},e_J\cdot{\bf v}_{{\mathcal I}})\cdot e_K\cdot{\bf v}_{{\mathcal I}}=0
$$
and hence proving \eqref{Jac22fore_8} reduces to proving that
\begin{equation}\label{Jac3fore_8}
\tilde L_2(e_J\cdot{\bf v}_{{\mathcal I}},e_K\cdot{\bf v}_{{\mathcal I}})\cdot e_I\cdot{\bf v}_{{\mathcal I}}+
\tilde L_2(e_K\cdot{\bf v}_{{\mathcal I}},e_I\cdot{\bf v}_{{\mathcal I}})\cdot e_J\cdot{\bf v}_{{\mathcal I}}
=0
\end{equation}
for all  subsets $I,J,K$ of $\{1,2,3,4,5,6,7,8\}$ satisfying one of (i), (ii), (iii) or (iv).
\vskip0,1cm
\noindent{\bf Case (i):} We have $I=J=\emptyset$ and $K=\{1,2,3,4,5,6,7,8\}$.
By Proposition \ref{explicitformforL} (iii) this means that
$$
\tilde L_2(e_J\cdot{\bf v}_{{\mathcal I}},e_K\cdot{\bf v}_{{\mathcal I}})\cdot e_I\cdot{\bf v}_{{\mathcal I}}=
-4B({\bf v}_{{\mathcal I}},e_K\cdot{\bf v}_{{\mathcal I}}){\bf v}_{{\mathcal I}}
$$
and
$$
\tilde L_2(e_K\cdot{\bf v}_{{\mathcal I}},{\bf v}_{{\mathcal I}})\cdot e_J\cdot{\bf v}_{{\mathcal I}}=
4B(e_K\cdot{\bf v}_{{\mathcal I}},{\bf v}_{{\mathcal I}}){\bf v}_{{\mathcal I}}.
$$
Hence
$$
\tilde L_2(e_J\cdot{\bf v}_{{\mathcal I}},e_K\cdot{\bf v}_{{\mathcal I}})\cdot e_I\cdot{\bf v}_{{\mathcal I}}+
\tilde L_2(e_K\cdot{\bf v}_{{\mathcal I}},e_I\cdot{\bf v}_{{\mathcal I}})\cdot e_J\cdot{\bf v}_{{\mathcal I}}
=0
$$
since $B$ is symmetric.
\vskip0,1cm
\noindent{\bf Case (ii):} We have $I=\emptyset$ and without loss of generality we can suppose
$J=\{1,2\}$ and $K=\{3,4,5,6,7,8\}$. By Proposition \ref{explicitformforL}(iii)  this means
$$
\tilde L_2(e_J\cdot{\bf v}_{{\mathcal I}},e_K\cdot{\bf v}_{{\mathcal I}})\cdot e_I\cdot{\bf v}_{{\mathcal I}}=
-2B(e_{12}{\bf v}_{{\mathcal I}},e_{345678}\cdot{\bf v}_{{\mathcal I}}){\bf v}_{{\mathcal I}}
$$
and by Proposition \ref{explicitformforL}(iv) that
$$
\tilde L_2(e_K\cdot{\bf v}_{{\mathcal I}},e_I\cdot{\bf v}_{{\mathcal I}})\cdot e_J\cdot {\bf v}_{{\mathcal I}}=
2B(e_{21}\cdot e_{345678}\cdot {\bf v}_{{\mathcal I}},{\bf v}_{{\mathcal I}}){\bf v}_{{\mathcal I}}.
$$
Hence from $B(e_{12}\cdot{\bf v}_{{\mathcal I}},e_{345678}\cdot{\bf v}_{{\mathcal I}})=B(e_{2}\cdot{\bf v}_{{\mathcal I}},e_1e_{345678}\cdot{\bf v}_{{\mathcal I}})=B({\bf v}_{{\mathcal I}},e_2e_1e_{345678}\cdot{\bf v}_{{\mathcal I}})$ and the fact that
$B$ is symmetric it follows that
$$
\tilde L_2(e_J\cdot{\bf v}_{{\mathcal I}},e_K\cdot{\bf v}_{{\mathcal I}})\cdot e_I\cdot{\bf v}_{{\mathcal I}}+
\tilde L_2(e_K\cdot{\bf v}_{{\mathcal I}},e_I\cdot{\bf v}_{{\mathcal I}})\cdot e_J\cdot{\bf v}_{{\mathcal I}}
=0
$$
\vskip0,1cm
\noindent{\bf Case (iii):} We have $I=\emptyset$ and without loss of generality we can suppose
$J=\{1,2,3,4\}$ and $K=\{5,6,7,8\}$. By Proposition \ref{explicitformforL} (iii) this means ($4=\frac{8}{2}!$) that
$$
\tilde L_2(e_J\cdot{\bf v}_{{\mathcal I}},e_K\cdot{\bf v}_{{\mathcal I}})\cdot e_I\cdot{\bf v}_{{\mathcal I}}=0
$$
and by Proposition \ref{explicitformforL} (ii) that
$$
\tilde L_2(e_K\cdot{\bf v}_{{\mathcal I}},{\bf v}_{{\mathcal I}})\cdot e_J\cdot{\bf v}_{{\mathcal I}}=0.
$$
It follows immediately that
$$
\tilde L_2(e_J\cdot{\bf v}_{{\mathcal I}},e_K\cdot{\bf v}_{{\mathcal I}})\cdot e_I\cdot{\bf v}_{{\mathcal I}}+
\tilde L_2(e_K\cdot{\bf v}_{{\mathcal I}},e_I\cdot{\bf v}_{{\mathcal I}})\cdot e_J\cdot{\bf v}_{{\mathcal I}}
=0
$$
\vskip0,1cm
\noindent{\bf Case (iv):} We can suppose without loss of generality that  $I=\{1,2\}$,
$J=\{3,4\}$ and $K=\{5,6,7,8\}$. By Proposition \ref{explicitformforL}(iv)  this means that
$$
\tilde L_2(e_J\cdot{\bf v}_{{\mathcal I}},e_K\cdot{\bf v}_{{\mathcal I}})\cdot  e_I\cdot{\bf v}_{{\mathcal I}}=
-2B(e_{12}e_{34}\cdot{\bf v}_{{\mathcal I}},e_{5678}\cdot{\bf v}_{{\mathcal I}}){\bf v}_{{\mathcal I}}
$$
and that
$$
\tilde L_2(e_K\cdot{\bf v}_{{\mathcal I}},e_I\cdot{\bf v}_{{\mathcal I}})\cdot e_J\cdot{\bf v}_{{\mathcal I}}=
2B(e_{43}e_{5678}\cdot{\bf v}_{{\mathcal I}},e_{12}\cdot{\bf v}_{{\mathcal I}}){\bf v}_{{\mathcal I}}.
$$
Hence from $B(e_{12}e_{34}\cdot{\bf v}_{{\mathcal I}},e_{5678}\cdot{\bf v}_{{\mathcal I}})=B({\bf v}_{{\mathcal I}},e_4e_3e_2e_1e_{5678}\cdot{\bf v}_{{\mathcal I}})$ and the fact that
$B$ is symmetric it follows that
$$
\tilde L_2(e_J\cdot{\bf v}_{{\mathcal I}},e_K\cdot{\bf v}_{{\mathcal I}})\cdot e_I\cdot{\bf v}_{{\mathcal I}}+
\tilde L_2(e_K\cdot{\bf v}_{{\mathcal I}},e_I\cdot{\bf v}_{{\mathcal I}})\cdot e_J\cdot{\bf v}_{{\mathcal I}}
=0.
$$
\begin{prop}
The Lie algebra constructed above is simple.
\end{prop}
\begin{demo} The  Lie bracket we have just defined on  the 248-dimensional vector space
$$
E=C^2(V,g)\oplus S_1
$$
 has the following properties:
\begin{itemize}
\item[(i)]
$C^2(V,g)$ is a simple Lie subalgebra;
\item[(ii)]
the bracket of $C^2(V,g)$ with $S_1$  defines a nontrivial  irreducible representation of $C^2(V,g)$ on $S_1$;
\item[(iii)]
${\rm dim}(C^2(V,g))>{\rm dim}(S_1)$;
\item[(iv)]
$[S_1,S_1]= C^2(V,g)$ .
\end{itemize}

We denote by $\pi_{C^2}:E\rightarrow C^2(V,g)$ and $\pi_{S_1}:E\rightarrow S_1$  the projections defined by the direct sum decomposition $E=C^2(V,g)\oplus S_1$. Clearly these are $C^2(V,g)$-equivariant maps whose respective kernels are $S_1$ and $C^2(V,g)$.

If ${\mathcal I}$ is a nonzero  ideal in $E$ then $\pi_{C^2}({\mathcal I})$ cannot be equal to $\{0\}$ - if so then  ${\mathcal I}\subseteq S_1$ which would imply ${\mathcal I}=S_1$ (cf. (ii)) and this  is impossible since $S_1$ is not an ideal (cf. (iv)). From the equivariance of $\pi_{C^2}$ and the irreducibility of $C^2(V,g)$ it then follows that 
\begin{equation}\label{projC}
\pi_{C^2}({\mathcal I})=C^2(V,g)
\end{equation}
and similarly, since $C^2(V,g)$ is not an ideal (cf. (ii)), we have
\begin{equation}\label{projS}
\pi_{S_1}({\mathcal I})=S_1.
\end{equation}
The  rank theorem for $\pi_{C^2}:{\mathcal I}\rightarrow C^2(V,g)$ and \eqref{projC} imply
$$
{\rm dim}({\mathcal I})\ge {\rm dim}(C^2(V,g))
$$
and by (iii) this means
$$
{\rm dim}({\mathcal I}\cap C^2(V,g))>0
$$
which  by (i)  implies
\begin{equation}\label{IintersectC}
{\mathcal I}\cap C^2(V,g)=C^2(V,g).
\end{equation}
It now follows from \eqref{projS}, \eqref{IintersectC} and the rank theorem for $\pi_{S_1}:{\mathcal I}\rightarrow S_1$  that
$$
{\rm dim}({\mathcal I})={\rm dim}({\mathcal I}\cap C^2(V,g))+{\rm dim}(\pi_{S_1}({\mathcal I}))={\rm dim}(C^2(V,g))+{\rm dim}(S_1)={\rm dim}(E)
$$
and hence 
$$
{\mathcal I}=E.
$$
\end{demo}
\subsection{Construction of split $e_7$}\hfill

Let $(V,g)$ be a twelve-dimensional vector space with a nondegenerate hyperbolic symmetric bilinear form $g$. Choose a 64-dimensional space of spinors $S=S_1\oplus S_2$ and  an isomorphism  $C(V,g)\cong  End( S)$. Since $n=6=2$ mod 4, $B:S\times S\rightarrow  k$ is even antisymmetric and $\tilde L_2:S\times S\rightarrow C^2(V,g)$  is even symmetric  (Proposition \ref{spinornormsymmetries} and Proposition \ref{Lsymmetries}). Hence $(S_1,B,C^2(V,g),L)$ is a
symplectic representation of the Lie algebra $C^2(V,g)$ possessing a natural equivariant symmetric bilinear form $\tilde L_2$ with values in $C^2(V,g)$. 

Following \cite{Fa} or \cite{GSSR} we can define a Lie bracket on the $133$-dimensional space
$$
E=C^2(V,g)\oplus sl(2,k)\oplus S_1\otimes k^2
$$
if $\tilde L_2$ (or a multiple of $\tilde L_2$) satisfies the equation
\begin{align}\label{CS1fore_7}
\nonumber
\tilde L_2(\psi_1,\psi_2)\cdot\psi_3-\tilde L_2(\psi_1,\psi_3)\cdot\psi_2=
-B(\psi_1,\psi_2)\psi_3+B(&\psi_1,\psi_3)\psi_2+2B(\psi_2,\psi_3)\psi_1\\
&\quad\forall \psi_1,\psi_2,\psi_3\in S_1.
\end{align}
Choose a polarisation $V={\mathcal I}\oplus {\mathcal E}$  of $V$ such that a  pure spinor ${\bf v}_{{\mathcal I}}$
corresponding to ${\mathcal I}$ is in $S_1$ and bases $i_1,\cdots ,i_n$ and $e_1,\cdots ,e_n$ of respectively ${\mathcal I}$ and ${\mathcal E}$ such that
\begin{align}\label{Cliffrel3}
\nonumber
&i_ai_b+i_bi_a=0\\
\nonumber
&e_ae_b+e_be_a=0\\
&i_ae_b+e_bi_a=\delta_{ab}.
\end{align}
Then 
$$
S_1={\rm Vect}\{e_I\cdot{\bf v}_{{\mathcal I}}:\quad\vert I\vert\text{ is even}\}
$$
and to prove the identity \eqref{CS1fore_7} it is sufficient to prove that
\begin{align}\label{CS2fore_7}
\nonumber
&\tilde L_2(e_I\cdot{\bf v}_{{\mathcal I}},e_J\cdot{\bf v}_{{\mathcal I}})\cdot e_K\cdot{\bf v}_{{\mathcal I}}-
\tilde L_2e_I\cdot{\bf v}_{{\mathcal I}},e_K\cdot{\bf v}_{{\mathcal I}})\cdot e_J\cdot{\bf v}_{{\mathcal I}}=\\
&-B(e_I\cdot{\bf v}_{{\mathcal I}},e_J\cdot{\bf v}_{{\mathcal I}})e_K\cdot{\bf v}_{{\mathcal I}}+
B(e_I\cdot{\bf v}_{{\mathcal I}},e_K\cdot{\bf v}_{{\mathcal I}})e_J\cdot{\bf v}_{{\mathcal I}}+
2B(e_J\cdot{\bf v}_{{\mathcal I}},e_K\cdot{\bf v}_{{\mathcal I}})e_I\cdot{\bf v}_{{\mathcal I}}
\end{align}
for all even oriented subsets $I,J,K$ of $\{1,2,3,4,5,6\}$. 

By the r\'esum\'e above,  the two terms on the LHS of this equation vanish unless $I\cap J\cap K=\emptyset$ and $I^c\cap J^c\cap K^c=\emptyset$ and by Proposition \ref{matrixB}  the three terms on the RHS vanish if $I\cap J\cap K\not=\emptyset$. Hence we need only prove \eqref{CS2fore_7} for even subsets of $\{1,2,3,4,5,6\}$ satisfying $I\cap J\cap K=\emptyset$ and $I^c\cap J^c\cap K^c=\emptyset$ . But then, again by the r\'esum\'e, by changing the polarisation if necessary, we can always assume that
$$
I\cap J=\emptyset,\quad K=I^c\cap J^c
$$
and then all terms on the LHS of  \eqref{CS2fore_7} vanish unless one of the sets $I,J,K$ has $0$ or $2$ elements. 
So  in fact to prove \eqref{CS2fore_7} for all even oriented subsets $I,J,K$ of $\{1,2,3,4,5,6\}$ it is sufficient  to prove that
\begin{align}
\nonumber
&\tilde L_2(e_I\cdot{\bf v}_{{\mathcal I}},e_J\cdot{\bf v}_{{\mathcal I}})\cdot e_K\cdot{\bf v}_{{\mathcal I}}-
\tilde L_2(e_I\cdot{\bf v}_{{\mathcal I}},e_K\cdot{\bf v}_{{\mathcal I}})\cdot e_J\cdot{\bf v}_{{\mathcal I}}=\\
&-B(e_I\cdot{\bf v}_{{\mathcal I}},e_J\cdot{\bf v}_{{\mathcal I}})e_K\cdot{\bf v}_{{\mathcal I}}+
B(e_I\cdot{\bf v}_{{\mathcal I}},e_K\cdot {\bf v}_{{\mathcal I}})e_J\cdot{\bf v}_{{\mathcal I}}+
2B(e_J\cdot{\bf v}_{{\mathcal I}},e_K\cdot{\bf v}_{{\mathcal I}})e_I\cdot{\bf v}_{{\mathcal I}}
\end{align}
for all  oriented subsets $I,J,K$ of $\{1,2,3,4,5,6\}$ such that
\begin{itemize}
\item  $\vert I\vert, \vert J\vert$ and $\vert K\vert$ are even;
\item $\vert I\vert+\vert J\vert+\vert K\vert=6$;
\item $I\cap J= J\cap K= K\cap I=\emptyset$;
\item One of $\vert I\vert, \vert J\vert$ or $\vert K\vert$ is equal to $0$ or $2$.
\end{itemize}
Up to permutations of $I,J,K$ the only possibilities are 
\begin{itemize}
\item[(i)]  $\vert I\vert=0$, $\vert J\vert=0$ and $\vert K\vert=6$;
\item[(ii)] $\vert I\vert=0$, $\vert J\vert=2$ and $\vert K\vert=4$;
\item[(iii)]  $\vert I\vert=2$, $\vert J\vert=2$ and $\vert K\vert=2$.
\end{itemize}
\vskip0,1cm
\noindent{\bf Case (i):} We have $I=J=\emptyset$ and  $K=\{1,2,3,4,5,6\}$. By Proposition \ref{explicitformforL}(ii)  this means that
$$
\tilde L_2(e_I\cdot{\bf v}_{{\mathcal I}},e_J\cdot{\bf v}_{{\mathcal I}})\cdot e_K\cdot{\bf v}_{{\mathcal I}}=0
$$
and by Proposition \ref{explicitformforL}(iii) that
$$
\tilde L_2(e_I\cdot{\bf v}_{{\mathcal I}},e_K\cdot{\bf v}_{{\mathcal I}})\cdot e_J\cdot{\bf v}_{{\mathcal I}}=
-3B({\bf v}_{{\mathcal I}},e_{123456}\cdot{\bf v}_{{\mathcal I}}){\bf v}_{{\mathcal I}}.
$$
Hence the LHS of \eqref{CS2fore_7} is 
$$
3B({\bf v}_{{\mathcal I}},e_{123456}\cdot{\bf v}_{{\mathcal I}}){\bf v}_{{\mathcal I}}
$$
which is equal to the RHS of \eqref{CS2fore_7} since by Proposition \ref{matrixB},
$$
-B(e_I\cdot{\bf v}_{{\mathcal I}},e_J\cdot{\bf v}_{{\mathcal I}})e_K\cdot{\bf v}_{{\mathcal I}}+
B(e_I\cdot{\bf v}_{{\mathcal I}},e_K\cdot{\bf v}_{{\mathcal I}})e_J\cdot{\bf v}_{{\mathcal I}}+
2B(e_J\cdot{\bf v}_{{\mathcal I}},e_K\cdot{\bf v}_{{\mathcal I}})e_I\cdot{\bf v}_{{\mathcal I}}
$$
reduces to
$$
0+
B({\bf v}_{{\mathcal I}},e_{123456}\cdot{\bf v}_{{\mathcal I}}){\bf v}_{{\mathcal I}}+
2B({\bf v}_{{\mathcal I}},e_{123456}\cdot{\bf v}_{{\mathcal I}})e_I\cdot{\bf v}_{{\mathcal I}}=
3B({\bf v}_{{\mathcal I}},e_{123456}\cdot{\bf v}_{{\mathcal I}}){\bf v}_{{\mathcal I}}.
$$
\vskip0,1cm
\noindent{\bf Case (ii):} We have $I=\emptyset$ and  without loss of generality we can suppose $J=\{1,2\}$ and $K=\{3,4,5,6\}$. By Proposition \ref{explicitformforL}(ii)  this means that
$$
\tilde L_2(e_I\cdot{\bf v}_{{\mathcal I}},e_J\cdot{\bf v}_{{\mathcal I}})\cdot e_K\cdot{\bf v}_{{\mathcal I}}=0
$$
and by Proposition \ref{explicitformforL}(iv) that
$$
\tilde L_2(e_I\cdot{\bf v}_{{\mathcal I}},e_K\cdot{\bf v}_{{\mathcal I}})\cdot e_J\cdot{\bf v}_{{\mathcal I}}=
2B(e_{21}{\bf v}_{{\mathcal I}},e_{3456}\cdot{\bf v}_{{\mathcal I}}){\bf v}_{{\mathcal I}}.
$$
Hence the LHS of \eqref{CS2fore_7} is 
$$
-2B(e_{21}{\bf v}_{{\mathcal I}},e_{3456}\cdot{\bf v}_{{\mathcal I}}){\bf v}_{{\mathcal I}}
$$
which is equal to the RHS of \eqref{CS2fore_7} since by Proposition \ref{matrixB},
$$
-B(e_I\cdot{\bf v}_{{\mathcal I}},e_J\cdot{\bf v}_{{\mathcal I}})e_K\cdot{\bf v}_{{\mathcal I}}+
B(e_I\cdot{\bf v}_{{\mathcal I}},e_K\cdot{\bf v}_{{\mathcal I}})e_J\cdot{\bf v}_{{\mathcal I}}+
2B(e_J\cdot{\bf v}_{{\mathcal I}},e_K\cdot{\bf v}_{{\mathcal I}})e_I\cdot{\bf v}_{{\mathcal I}}
$$
reduces to
$$
0+0+2B(e_{12}.{\bf v}_{{\mathcal I}},e_{3456}\cdot{\bf v}_{{\mathcal I}}){\bf v}_{{\mathcal I}}=
-2B(e_{21}{\bf v}_{{\mathcal I}},e_{3456}\cdot{\bf v}_{{\mathcal I}}){\bf v}_{{\mathcal I}}.
$$
\vskip0,1cm
\noindent{\bf Case (iii):} Without loss of generality we can suppose $I=\{1,2\}$,  $J=\{3,4\}$ and $K=\{5,6\}$. By Proposition \ref{explicitformforL}(iv)  this means that
$$
\tilde L_2(e_I\cdot{\bf v}_{{\mathcal I}},e_J\cdot{\bf v}_{{\mathcal I}})\cdot e_K\cdot{\bf v}_{{\mathcal I}}=
2B(e_{65}e_{12}\cdot{\bf v}_{{\mathcal I}},e_{34}\cdot{\bf v}_{{\mathcal I}}){\bf v}_{{\mathcal I}}
$$
and  that
$$
\tilde L_2(e_I\cdot{\bf v}_{{\mathcal I}},e_K\cdot{\bf v}_{{\mathcal I}})\cdot e_J\cdot{\bf v}_{{\mathcal I}}=
2B(e_{43}e_{12}\cdot{\bf v}_{{\mathcal I}},e_{56}\cdot{\bf v}_{{\mathcal I}}){\bf v}_{{\mathcal I}}.
$$
Hence the LHS of \eqref{CS2fore_7} is 
$$
2B(e_{65}e_{12}\cdot{\bf v}_{{\mathcal I}},e_{34}\cdot{\bf v}_{{\mathcal I}}){\bf v}_{{\mathcal I}}-
2B(e_{43}e_{12}\cdot{\bf v}_{{\mathcal I}},e_{56}\cdot{\bf v}_{{\mathcal I}}){\bf v}_{{\mathcal I}}=0.
$$
The RHS of \eqref{CS2fore_7} also vanishes since each term of
$$
-B(e_I\cdot{\bf v}_{{\mathcal I}},e_J\cdot{\bf v}_{{\mathcal I}})e_K\cdot{\bf v}_{{\mathcal I}}+
B(e_I\cdot{\bf v}_{{\mathcal I}},e_K\cdot{\bf v}_{{\mathcal I}})e_J\cdot{\bf v}_{{\mathcal I}}+
2B(e_J\cdot{\bf v}_{{\mathcal I}},e_K\cdot{\bf v}_{{\mathcal I}})e_I\cdot{\bf v}_{{\mathcal I}}
$$
vanishes by Proposition \ref{matrixB}.
\begin{prop}
The Lie algebra constructed above is simple.
\end{prop}
\begin{demo}
The Lie bracket  we have just defined on the $133$-dimensional space
$$
E=C^2(V,g)\oplus sl(2,k)\oplus S_1\otimes k^2
$$
has the properties:
\begin{itemize}
\item[(i)]
$C^2(V,g)$ and $sl(2,k)$ are commuting simple Lie subalgebras;
\item[(ii)]
the bracket of $C^2(V,g)\oplus sl(2,k)$ with $S_1\otimes k^2$ defines a faithful, irreducible representation of $C^2(V,g)\oplus sl(2,k)$ on $S_1\otimes k^2$;
\item[(iii)]
$[S_1\otimes k^2,S_1\otimes k^2]= C^2(V,g)\oplus sl(2,k)$.
\end{itemize}
Pick any standard semisimple $h$ in $sl(2,k)$. Then ${\rm ad}(h):E\rightarrow E$ is diagonalisable with  eigenvalues $\{0,\pm1,\pm2\}$ and:
\begin{itemize}
\item
 $C^2(V,g)\oplus sl(2,k)=E_+$ is the sum of the eigenspaces corresponding to even eigenvalues;
 \item
  $S_1\otimes k^2=E_-$ is the sum of the eigenspaces corresponding to odd eigenvalues.
\end{itemize}
Let ${\mathcal I}$ be a nonzero ideal in $E$. Then $[h,{\mathcal I}]\subseteq{\mathcal I}$ and hence
$$
{\mathcal I}={\mathcal I}\cap E_+\oplus {\mathcal I}\cap E_-.
$$
If $ {\mathcal I}\cap E_-=\{0\}$ then ${\mathcal I}={\mathcal I}\cap E_+\not=\{0\}$ which is impossible since no nontrivial ideal of $E$ can be contained in $E_+$ ($E_+$ acts faithfully on $E_-$ by (ii)). Hence $ {\mathcal I}\cap E_-\not=\{0\}$.

If  $ {\mathcal I}\cap E_-\not=\{0\}$ then in fact  $ {\mathcal I}\cap E_-=E_-$ since ${\mathcal I}\cap E_-$ is stable under $E_+$ and  $E_-$ is an irreducible representation of $E_+$ (cf (ii)).
However if ${\mathcal I}$ contains $E_-$ it contains $E_+$ by (iii) and hence ${\mathcal I}=E$.

\end{demo}

\subsection{Construction of split $e_6$}\hfill

Let $(V,g)$ be a ten-dimensional vector space with a nondegenerate hyperbolic symmetric bilinear form $g$. Choose a 32-dimensional space of spinors $S=S_1\oplus S_2$ and  an isomorphism  $C(V,g)\cong End( S)$. Since $n=5=1$ (mod 4), $B:S\times S\rightarrow k$ is odd symmetric (Proposition \ref{spinornormsymmetries}), $\tilde L_2:S\times S\rightarrow C^2(V,g)$  is odd antisymmetric  (Proposition \ref{Lsymmetries}) and  $ L_{10}:S\times S\rightarrow C^{10}(V,g)$ is  odd antisymmetric (Proposition \ref{ L_{2n}symmetries}).

Consider the $78$-dimensional vector space
$$
E= C^{10}(V,g)\oplus C^2(V,g)\oplus S.
$$
Following the procedure in \S 7 define  $[\phantom x ,\phantom y ]:E\otimes E\rightarrow E$ to be the unique antisymmetric bilinear map such that:
\begin{align}
[A,B]&=AB-BA\quad&\hbox{ if }A,B\in C^2\oplus C^{10}\\
[A,\psi]&=A\cdot\psi\quad&\hbox{ if }A\in C^2\oplus C^{10},\,\psi\in S\\
[\psi_1,\psi_2]&=2\tilde L_2(\psi_1,\psi_2)+96 L_{10}(\psi_1,\psi_2)&\quad\hbox{ if }\psi_1,\psi_2\in S.
\end{align}
Since $C^2,C^2\oplus C^{10}$ are Lie algebras , since $S$
is a representation and  since $\tilde L_2:S\times S\rightarrow C^2$ and $2\tilde L_2+ 96L_{10}:S\times S\rightarrow C^2\oplus C^{10}$ are equivariant maps, this map defines a Lie bracket on $E$ if and only if the following Jacobi identies are satisfied:
\begin{align}
\nonumber
(2\tilde L_2+96L_{10})(\psi_1,\psi_2)\cdot\psi_3+(2\tilde L_2+96 L_{10})(\psi_2,\psi_3)\cdot\psi_1+(2\tilde L_2+96 L_{10})(\psi_3,\psi_1)\cdot\psi_2=0\\
\forall\psi_1,\psi_2,\psi_3\in S.
\end{align}
Recall that if  $\varepsilon\in C^{10}(V,g)$ is a grading operator then
$$
 L_{10}(\psi_1,\psi_2)=\frac{1}{2^5}B(\psi_1,\varepsilon\cdot\psi_2)\varepsilon
$$
so that this is equivalent to
\begin{align}\label{Jac1fore_6}
\nonumber
2\tilde L_2(\psi_1,\psi_2)\cdot\psi_3+2\tilde L_2(\psi_2,\psi_3)\cdot\psi_1+2\tilde L_2(\psi_3,\psi_1)\cdot\psi_2=&\\
\nonumber
-3B(\psi_1,\varepsilon\cdot\psi_2)\varepsilon\cdot\psi_3-
3B(\psi_2,\varepsilon\cdot\psi_3)\varepsilon\cdot\psi_1&-3B(\psi_3,\varepsilon\cdot\psi_1)\varepsilon\cdot\psi_2\\
&\forall\psi_1,\psi_2,\psi_3\in S.
\end{align}

A purist would say that in this example we should use the graded spinor norm described in \S 3.1. However, the calculations in \S 6 were  all done with the usual spinor norm. For the sake of ease of verifying the result we shall not use the graded norm but just the usual spinor norm.

Choose a polarisation $V={\mathcal I}\oplus {\mathcal E}$  of $V$ such that a  pure spinor ${\bf v}_{{\mathcal I}}$
corresponding to ${\mathcal I}$ is in $S_1$ and bases $i_1,\cdots ,i_n$ and $e_1,\cdots ,e_n$ of respectively ${\mathcal I}$ and ${\mathcal E}$ such that
\begin{align}
\nonumber
&i_ai_b+i_bi_a=0\\
\nonumber
&e_ae_b+e_be_a=0\\
&i_ae_b+e_bi_a=\delta_{ab}.
\end{align}
Then 
$$
S_1={\rm Vect}\{e_I\cdot{\bf v}_{{\mathcal I}}: \ \vert I\vert\text{ is even}\},\quad
S_2={\rm Vect}\{e_I\cdot{\bf v}_{{\mathcal I}}: \ \vert I\vert\text{ is odd}\}
$$
 and to prove the identity \eqref{Jac1fore_6} it is sufficient to prove that
\begin{align}\label{Jac2e_6}
\nonumber
&2\tilde L_2(e_I\cdot{\bf v}_{{\mathcal I}},e_J\cdot{\bf v}_{{\mathcal I}})\cdot e_K\cdot{\bf v}_{{\mathcal I}}+2\tilde L_2(e_J\cdot{\bf v}_{{\mathcal I}},e_K\cdot{\bf v}_{{\mathcal I}}).e_I\cdot{\bf v}_{{\mathcal I}}+2\tilde L_2(e_K\cdot{\bf v}_{{\mathcal I}},e_I\cdot{\bf v}_{{\mathcal I}})\cdot e_J\cdot{\bf v}_{{\mathcal I}}=\\
&-3B(e_I\cdot{\bf v}_{{\mathcal I}},\varepsilon\cdot e_J\cdot{\bf v}_{{\mathcal I}})\varepsilon\cdot e_K\cdot{\bf v}_{{\mathcal I}}-
3B(e_J\cdot{\bf v}_{{\mathcal I}},\varepsilon\cdot e_K\cdot{\bf v}_{{\mathcal I}})\varepsilon\cdot e_I\cdot{\bf v}_{{\mathcal I}}-3B(e_K\cdot{\bf v}_{{\mathcal I}},\varepsilon\cdot  e_I\cdot{\bf v}_{{\mathcal I}})\varepsilon\cdot e_J\cdot{\bf v}_{{\mathcal I}}
\end{align}
for all oriented subsets $I,J,K$ of $\{1,2,3,4,5\}$.

If $\vert I\vert,\vert J\vert ,\vert K\vert $ are of the same parity, all terms in this expression vanish ($\tilde L_2$ and $B$ are odd) and the identity is true. If $\vert I\vert,\vert J\vert ,\vert K\vert $ are not of the same parity two of them must be even and one odd, and without loss of generality we can suppose $\vert I\vert,\vert J\vert$ are even and $\vert K\vert $ is odd.  Since $\tilde L_2$ and $B$ are odd the identity \eqref{Jac2e_6} then reduces to
\begin{align}\label{Jac3e_6}
\nonumber
2\tilde L_2(e_J\cdot{\bf v}_{{\mathcal I}},e_K\cdot{\bf v}_{{\mathcal I}})\cdot &e_I\cdot {\bf v}_{{\mathcal I}}+2\tilde L_2(e_K\cdot{\bf v}_{{\mathcal I}},e_I\cdot{\bf v}_{{\mathcal I}})\cdot e_J\cdot{\bf v}_{{\mathcal I}}=\\
3B(e_J\cdot{\bf v}_{{\mathcal I}},&\varepsilon\cdot e_K\cdot{\bf v}_{{\mathcal I}})\varepsilon\cdot e_I\cdot{\bf v}_{{\mathcal I}}-3B(e_K\cdot{\bf v}_{{\mathcal I}},\varepsilon\cdot e_I\cdot{\bf v}_{{\mathcal I}})\varepsilon\cdot e_J\cdot{\bf v}_{{\mathcal I}}.
\end{align}
By the r\'esum\'e above,  the two terms on the LHS of this equation vanish unless $I\cap J\cap K=\emptyset$ and $I^c\cap J^c\cap K^c=\emptyset$ and by Proposition \ref{matrixB}  the three terms on the RHS vanish if $I\cap J\cap K\not=\emptyset$. 

If $I\cap J\cap K=\emptyset$ and $I^c\cap J^c\cap K^c=\emptyset$ then by changing the polarisation if necessary, we can always assume that
$$
I\cap J=\emptyset,\quad K=I^c\cap J^c,
$$
and then all terms on the LHS of  \eqref{Jac3e_6} vanish unless one of the sets $I,J,K$ has $0$ or $2$ elements. 
Up to permutations of $I,J,K$ the only possibilities are 
\begin{itemize}
\item[(i)]  $\vert I\vert=0$, $\vert J\vert=0$ and $\vert K\vert=5$;
\item[(ii)] $\vert I\vert=0$, $\vert J\vert=2$ and $\vert K\vert=3$;
\item[(iii)]  $\vert I\vert=0$, $\vert J\vert=4$ and $\vert K\vert=1$;
\item[(iv)]  $\vert I\vert=2$, $\vert J\vert=2$ and $\vert K\vert=1$.
\end{itemize}
\vskip0,1cm
\noindent{\bf Case (i):} We have $I=J=\emptyset$ and  $K=\{1,2,3,4,5\}$. By Proposition \ref{explicitformforL}(iii)  this means that
$$
2\tilde L_2(e_J\cdot{\bf v}_{{\mathcal I}},e_K\cdot{\bf v}_{{\mathcal I}})\cdot e_I\cdot{\bf v}_{{\mathcal I}}=
2B({\bf v}_{{\mathcal I}},e_{12345}\cdot{\bf v}_{{\mathcal I}})(0-\frac{5}{2}){\bf v}_{{\mathcal I}}
$$
and  that
$$
2\tilde L_2(e_K\cdot{\bf v}_{{\mathcal I}},e_I\cdot{\bf v}_{{\mathcal I}})\cdot e_J\cdot{\bf v}_{{\mathcal I}}=
2B(e_{12345}\cdot{\bf v}_{{\mathcal I}},{\bf v}_{{\mathcal I}})(5-\frac{5}{2}){\bf v}_{{\mathcal I}}
$$
Hence the LHS of \eqref{Jac3e_6} vanishes since $B$ is symmetric as does the RHS for the same reason.
\vskip0,1cm
\noindent{\bf Case (ii):} We have $I=\emptyset$ and  without loss of generality we can suppose $J=\{1,2\}$ and $K=\{3,4,5\}$. By Proposition \ref{explicitformforL}(iii)  this means that
$$
2\tilde L_2(e_J\cdot{\bf v}_{{\mathcal I}},e_K\cdot{\bf v}_{{\mathcal I}})\cdot e_I\cdot{\bf v}_{{\mathcal I}}=
2B(e_{12}{\bf v}_{{\mathcal I}},e_{345}\cdot{\bf v}_{{\mathcal I}})(2-\frac{5}{2}){\bf v}_{{\mathcal I}},
$$
and  by Proposition \ref{explicitformforL}(iv)  that
$$
2\tilde L_2(e_K\cdot{\bf v}_{{\mathcal I}},e_I\cdot{\bf v}_{{\mathcal I}})\cdot e_J\cdot{\bf v}_{{\mathcal I}}=
4B(e_{21}\cdot e_{345}\cdot{\bf v}_{{\mathcal I}},{\bf v}_{{\mathcal I}}){\bf v}_{{\mathcal I}}.
$$
Hence the LHS of \eqref{Jac3e_6} reduces to
$$
3B(e_{21}\cdot e_{345}\cdot{\bf v}_{{\mathcal I}},{\bf v}_{{\mathcal I}}){\bf v}_{{\mathcal I}}.
$$
The RHS of \eqref{Jac3e_6} is
$$
3B(e_J\cdot{\bf v}_{{\mathcal I}},\varepsilon\cdot e_K\cdot{\bf v}_{{\mathcal I}})\varepsilon\cdot e_I\cdot{\bf v}_{{\mathcal I}}-3B(e_K\cdot{\bf v}_{{\mathcal I}},\varepsilon\cdot e_I\cdot{\bf v}_{{\mathcal I}})\varepsilon\cdot e_J\cdot{\bf v}_{{\mathcal I}}
$$
which by Proposition \ref{matrixB} also reduces to
$$
3B(e_{12}\cdot{\bf v}_{{\mathcal I}},e_{345}\cdot {\bf v}_{{\mathcal I}}){\bf v}_{{\mathcal I}}-0=3B(e_{21}\cdot e_{345}\cdot {\bf v}_{{\mathcal I}},{\bf v}_{{\mathcal I}}){\bf v}_{{\mathcal I}}.
$$
\vskip0,1cm
\noindent{\bf Case (iii):} We have $I=\emptyset$ and  without loss of generality we can suppose $J=\{1,2,3,4\}$ and $K=\{5\}$. By Proposition \ref{explicitformforL}(iii)  this means that
$$
2\tilde L_2(e_J\cdot{\bf v}_{{\mathcal I}},e_K\cdot{\bf v}_{{\mathcal I}})\cdot e_I\cdot{\bf v}_{{\mathcal I}}=
2B(e_{1234}\cdot{\bf v}_{{\mathcal I}},e_{5}\cdot{\bf v}_{{\mathcal I}})(4-\frac{5}{2}){\bf v}_{{\mathcal I}}
$$
and by Proposition \ref{explicitformforL}(ii)  that
$$
2\tilde L_2(e_K\cdot{\bf v}_{{\mathcal I}},e_I\cdot{\bf v}_{{\mathcal I}})\cdot e_J\cdot{\bf v}_{{\mathcal I}}=0
$$
Hence the LHS of \eqref{Jac3e_6} reduces to
$$
3B(e_{12345}\cdot{\bf v}_{{\mathcal I}},{\bf v}_{{\mathcal I}}){\bf v}_{{\mathcal I}}.
$$
The RHS of \eqref{Jac3e_6} is
$$
3B(e_J\cdot{\bf v}_{{\mathcal I}},\varepsilon\cdot e_K\cdot{\bf v}_{{\mathcal I}})\varepsilon\cdot e_I\cdot{\bf v}_{{\mathcal I}}-3B(e_K\cdot{\bf v}_{{\mathcal I}},\varepsilon\cdot e_I\cdot{\bf v}_{{\mathcal I}})\varepsilon\cdot e_J\cdot{\bf v}_{{\mathcal I}}
$$
which by Proposition \ref{matrixB} also reduces to
$$
3B(e_{1234}{\bf v}_{{\mathcal I}},e_{5}\cdot{\bf v}_{{\mathcal I}}){\bf v}_{{\mathcal I}}-0=3B(e_{12345}\cdot{\bf v}_{{\mathcal I}},{\bf v}_{{\mathcal I}}){\bf v}_{{\mathcal I}}.
$$
\vskip0,1cm
\noindent{\bf Case (iv):} We can suppose  without loss of generality that $I=\{1,2\}$, $J=\{3,4\}$ and $K=\{5\}$. By Proposition \ref{explicitformforL}(iv)  this means 
$$
2\tilde L_2(e_J\cdot{\bf v}_{{\mathcal I}},e_K\cdot{\bf v}_{{\mathcal I}})\cdot e_I\cdot{\bf v}_{{\mathcal I}}=
4B(e_{21}\cdot e_{34}\cdot {\bf v}_{{\mathcal I}},e_{5}\cdot{\bf v}_{{\mathcal I}}){\bf v}_{{\mathcal I}}
$$
and 
$$
2\tilde L_2(e_K\cdot{\bf v}_{{\mathcal I}},e_I\cdot{\bf v}_{{\mathcal I}})\cdot e_J\cdot{\bf v}_{{\mathcal I}}=
4B(e_{43}\cdot e_{5}\cdot{\bf v}_{{\mathcal I}},e_{12}\cdot{\bf v}_{{\mathcal I}}){\bf v}_{{\mathcal I}}.
$$
Since
$$
B(e_{21}\cdot e_{34}\cdot{\bf v}_{{\mathcal I}},e_{5}\cdot{\bf v}_{{\mathcal I}}){\bf v}_{{\mathcal I}}=
B(e_{34}\cdot e_{21}\cdot {\bf v}_{{\mathcal I}},e_{5}\cdot{\bf v}_{{\mathcal I}}){\bf v}_{{\mathcal I}}=
B(e_{21}\cdot{\bf v}_{{\mathcal I}},e_{43}\cdot e_5 \cdot{\bf v}_{{\mathcal I}}){\bf v}_{{\mathcal I}}
$$
and $B$ is symmetric, the LHS of \eqref{Jac3e_6}  vanishes as do
all terms on the RHS by Proposition \ref{matrixB}.
\begin{prop}
The Lie algebra constructed above is simple.
\end{prop}
\begin{demo}
The  Lie bracket we have just defined on the $78$-dimensional space
$$
E=C^2(V,g)\oplus C^{10}(V,g)\oplus S
$$
has the properties:
\begin{itemize}
\item[(i)]
$C^2(V,g)$ is a simple Lie subalgebra and  $C^{10}(V,g)$ is its one-dimensional commutant in $E$;
\item[(ii)]
there exists $\varepsilon\in C^{10}(V,g)$ such that;
${\rm ad}(\varepsilon):E\rightarrow E$ is diagonalisable with  eigenvalues $\{0,\pm1\}$ and $S=E_{-1}\oplus E_1$ is the decomposition
of $S$ into faithful, non-isomorphic irreducible $C^2(V,g)$-modules.
\item[(iii)]
$[E_{-1}, E_1]= C^2(V,g)\oplus C^{10}(V,g)$.
\end{itemize}
Let ${\mathcal I}$ be a nonzero ideal in $E$. Then $[\varepsilon,{\mathcal I}]\subseteq{\mathcal I}$ and hence
$$
{\mathcal I}={\mathcal I}\cap E_0\oplus {\mathcal I}\cap E_{-1}\oplus  {\mathcal I}\cap E_1.
$$
If $ {\mathcal I}\cap E_{-1}= {\mathcal I}\cap E_{1}=\{0\}$ then ${\mathcal I}={\mathcal I}\cap E_0\not=\{0\}$ which is impossible since no nontrivial ideal of $E$ can be contained in $E_0$ by (ii) ($E_0$ acts faithfully on $S$ by (ii)). Hence either $ {\mathcal I}\cap E_{-1}\not=\{0\}$ or $ {\mathcal I}\cap E_{1}\not=\{0\}$.

If  $ {\mathcal I}\cap E_{-1}\not=\{0\}$ then in fact  $ {\mathcal I}\cap E_{-1}=E_{-1}$ since ${\mathcal I}\cap E_{-1}$ is stable under $E_0$ and  $E_{-1}$ is an irreducible representation of $E_0$ (cf (ii)).
However if ${\mathcal I}$ contains $E_{-1}$ it contains $E_0$ by (iii) and hence ${\mathcal I}=E$. Similarly, if $ {\mathcal I}\cap E_{1}\not=\{0\}$ then 
${\mathcal I}=E$ and the proposition is proved.
\end{demo}

\end{document}